%% file: Hewett_FLP_v1_preprint.tex
\documentclass[preprint,12pt]{elsarticle}
\journal{.}

\usepackage{amsmath}
\usepackage{amssymb}
\usepackage{enumerate}
\usepackage{cases}
\usepackage{multirow}
\usepackage{subfigure}
\usepackage{graphicx}
\usepackage{appendix}
\usepackage{pgf,tikz}
\usetikzlibrary{arrows}

\newif\ifpre
\pretrue
\begin{document}
\newcommand{\hx}{\hat{x}}
\newcommand{\hy}{\hat{y}}
\newcommand{\cy}{\check{y}}
\newcommand{\ty}{\tilde{y}}
\newcommand{\tg}{\tilde{g}}
\newcommand{\hp}{\hat{p}}
\newcommand{\hP}{\hat{P}}
\newcommand{\hq}{\hat{q}}
\newcommand{\hQ}{\hat{Q}}
\newcommand{\hn}{\hat{n}}
\newcommand{\hmu}{\hat{\mu}}
\newcommand{\tsigma}{\tilde{\sigma}}
\newcommand{\tSigma}{\tilde{\Sigma}}
\newcommand{\hV}{\hat{V}}
\newcommand{\rA}{{\rm A}}
\newcommand{\Fr}{{\rm Fr}}
\newcommand{\Heaviside}{H}
\newcommand{\Ai}{{\rm Ai}}
\newcommand{\tf}[1]{{\begin{scriptsize} \it #1 \end{scriptsize}}}
\newcommand{\Bi}{{\rm Bi}}
\newcommand{\F}{{Fig.~}}
\newcommand{\Fs}{{Figs~}}
\newcommand{\cF}{\mathcal{F}}
\newcommand{\tagR}[1]{\tag{R.\ref{#1}}}
\newcommand{\tagN}[1]{\tag{N.\ref{#1}}}
\newcommand{\includegraphicsDave}[2][]{\includegraphics[#1]{#2_LOWRES}}		
\newcommand{\mmbox}[1]{\fbox{\ensuremath{\displaystyle{ #1 }}}}	
\input{macros}
\begin{frontmatter}
\title{Tangent ray diffraction and the Pekeris caret function}

\author{D.\ P.\ Hewett}
\address{Mathematical Institute, University of Oxford, Radcliffe Observatory Quarter, Woodstock Road, Oxford, 
OX2 6GG, UK, {\tt hewett\char'100maths.ox.ac.uk}, +44 (0)1865 270744.}

\begin{abstract}

We study the classical problem of high frequency scattering of an incident plane wave by a smooth convex two-dimensional body. We present a new integral representation of the leading order solution in the ``Fock region'', i.e. the neighbourhood of a point of tangency between the incident rays and the scatterer boundary, from which the penumbra (light-shadow boundary) effects originate. 
The new representation, which is equivalent to the classical Fourier integral representation and its well-studied ``forked contour'' regularisation, reveals that the Pekeris caret function (sometimes referred to as a ``Fock-type integral'' or a ``Fock scattering function''), a special function already known to describe the field in the penumbra, is also an intrinsic part of the solution in the inner Fock region. 
We also provide the correct interpretation of a divergent integral arising in the analysis of Tew et al.\ (Wave Motion 32, 2000), enabling the results of that paper to be used for quantitative calculations.

\end{abstract}
\end{frontmatter}
%

\section{Introduction}
\label{sec:intro}

High frequency time-harmonic scattering of an incident wave by a smooth convex body is a classical problem in linear wave propagation, and 
has been the subject of intense study over the past 70 or so years. The earliest published mathematical studies seem to be those of Fock \cite{Fo:46}, Fock and Leontovich \cite{LeFo:45}, and Pekeris \cite{Pe:47}; generalisations and reinterpretations of these early works have subsequently been given by numerous different authors: for a comprehensive review of the literature see e.g.\ \cite{Fo:65,BaKi:79,BaBu:91,TeChKiOcSmZa:00} and the many references therein. Other notable works include \cite{Br:66,Lu:67,He:68,MeTa:86,BuLy:87}. 

\begin{figure}[t!]
\begin{center}
\includegraphics[scale=0.5]{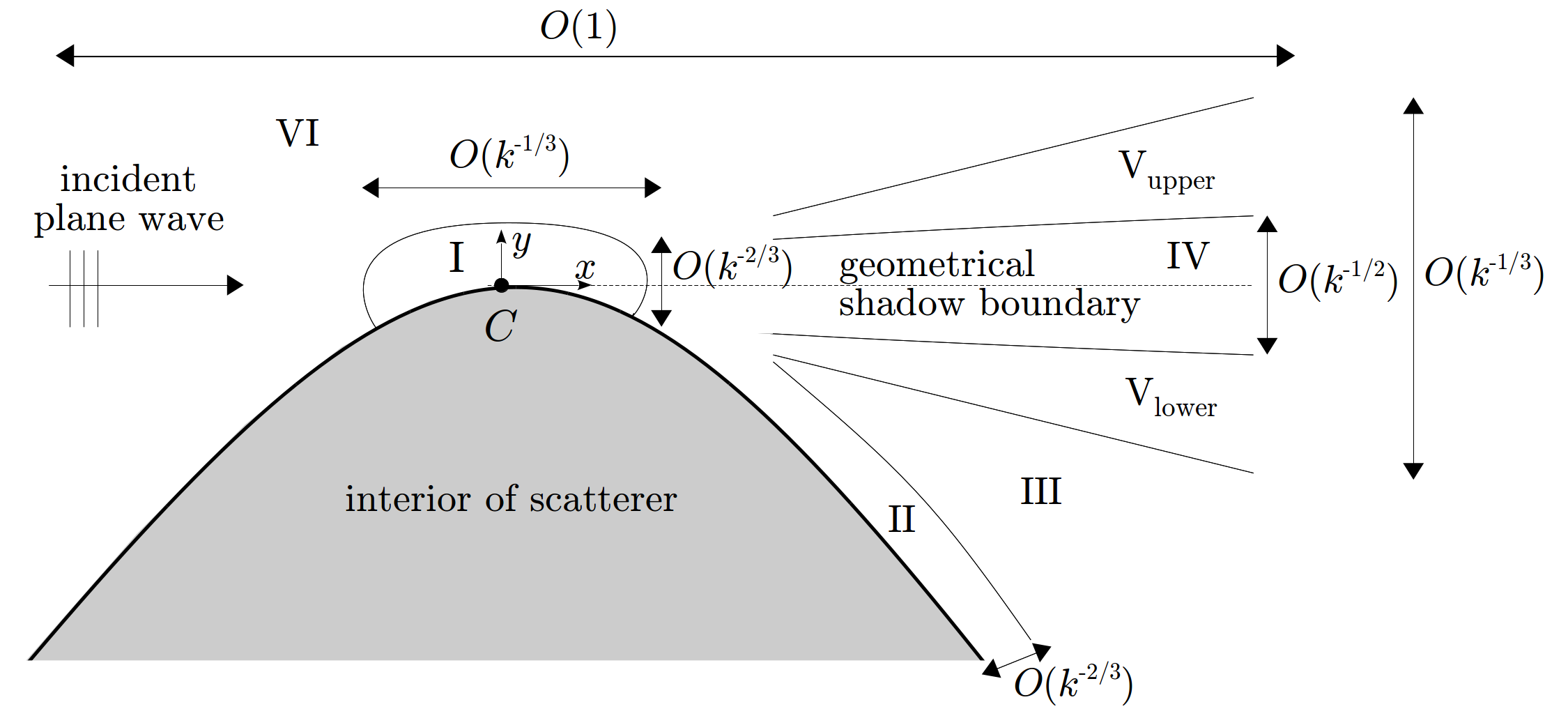}
\caption{
High frequency scattering by a smooth convex body: schematic (based on \cite[Fig.~2]{TeChKiOcSmZa:00}) showing the asymptotic structure of the field near a point of ray tangency (labelled $C$) and the geometrical shadow boundary. 
Here $k$ denotes the (nondimensional) wavenumber (cf.\ \rf{eqn:HE}). 
The inner Fock region (labelled I) is a small neighbourhood of the tangency point of width $\ord{k^{-1/3}}$ and height $\ord{k^{-2/3}}$. The ``Airy layer" or creeping wave region (labelled II) is a thin layer of height $\ord{k^{-2/3}}$ in the deep shadow close to the boundary, in which the field is described by creeping waves, which propagate along the boundary away from the point $C$. The amplitude of these creeping waves decays exponentially as a function of distance around the boundary from $C$, and the creeping waves shed rays tangentially into the deep shadow region away from the boundary (labelled III). 
The transition between this deep shadow region and the illuminated region (labelled VI), in which the specularly reflected field is present, takes place in the penumbra, which has width $\ord{k^{-1/3}}$, and which itself has a three-layer sub-structure. Around the geometrical shadow boundary (shown as a dashed line), there is an inner parabolic region of width $\ord{k^{-1/2}}$ (labelled IV) in which the field is described by a Fresnel integral. Here the field is (to leading order) independent of the boundary condition imposed on the scatterer boundary. Above and below this Fresnel region there are transition regions (labelled V$_\textrm{upper}$ and V$_\textrm{lower}$) across which the Fresnel integral solution matches respectively to the solutions in the shed creeping and illuminated regions.
}
\label{fig:OuterRegions}
\end{center}
\end{figure}

It is well known that, for the classical two-dimensional problem with plane wave incidence, the propagation domain around a point of tangency between the incident rays and the scatterer boundary can be divided into different regions, within each of which the leading order high frequency approximation of the wave field takes different forms. A concise and systematic treatment of this problem in the framework of matched asymptotic expansions can be found in \cite{TeChKiOcSmZa:00}; 
a schematic illustrating the structure of the resulting wave field 
is shown in \F\ref{fig:OuterRegions}. 
%

The purpose of the current article is not to provide yet another review of this well-known material, but rather to report a new representation for the solution in the inner ``Fock region'' (labelled I in \F\ref{fig:OuterRegions}), which provides a fresh perspective on this classical problem. 
As was first noted by Fock and Leontovich in the 1940s (see e.g.\ \cite{LeFo:45,Fo:46}), the scattered field close to a tangency point can be approximated by the product of the incident plane wave with a certain solution of the parabolic wave equation (PWE). 
In Fock and Leontovich's original work, an exact formula for the PWE solution was given in the form of a Fourier integral involving products and ratios of Airy functions (see \rf{eqn:FLScatClassical} below). The Fourier integral is defined improperly and does not converge absolutely, but Babich and coworkers subsequently showed how it could be regularised by deforming the integration contour into the complex plane (see e.g.\ \cite{BaKi:79,BaBu:91}). Unfortunately, a standard contour deformation is not sufficient; instead one has to decompose the integrand into two terms and deform the integration contour differently for each, leading to a so-called ``forked contour'' regularisation (see \rf{eqn:FLScatHorny} below).

Our new solution representation, presented in \rf{eqn:FLScatNew} below, offers an alternative regularisation of the classical Fourier integral solution. It takes the form of a single rapidly convergent contour integral, whose integrand comprises an exponential factor, with polynomial exponent, multiplied by a special function which, following \cite{Lo:59,James}, we refer to as the ``Pekeris caret function''. 
As will be explained in \S\ref{sec:NewSoln}, this function is closely related 
to an integral apparently first introduced by Pekeris \cite{Pe:47} in the 1940s. However, since the same integral also appeared around the same time in the work of Fock \cite[Eqn (4.03)]{Fo:65}, what we call the Pekeris caret function is referred to by some as a ``Fock-type integral'' \cite{Pe:87} or ``Fock scattering function'' \cite{Sa:09}. It is already well known (see \cite{Pe:47}, \cite[\S4]{Fo:65}, \cite[\S6.8]{BaKi:79}, \cite[\S13.7]{BaBu:91} and \cite[\S2.5]{TeChKiOcSmZa:00}) that this special function describes the field in the transition regions of the penumbra (regions V$_\textrm{upper}$ and V$_\textrm{lower}$ of \F\ref{fig:OuterRegions}). 
What our new representation reveals is that the Pekeris caret function is also an intrinsic part of the solution in the inner Fock region, something which is not at all obvious from the classical solution representation and its forked contour regularisation. As such, we argue that the Pekeris caret function should be regarded as the key special function describing tangent ray diffraction. 
As we shall explain, the specularly reflected field, the penumbra fields, and the creeping field in the deep shadow all emerge from our new contour integral representation by straightforward applications of the steepest descent (saddle point) method, combined with knowledge of certain elementary properties of the Pekeris caret function. 
In particular, the smooth switching on/off of the incident wave as one moves across the penumbra manifests itself in the coalescence of a saddle point with the single simple pole of the Pekeris caret function.

A secondary objective of the paper is to provide an interpretation (in terms of the Pekeris caret function) of a divergent integral appearing in \cite[Eqn (2.63)]{TeChKiOcSmZa:00} in the formula for the behaviour in the transition regions V$_\textrm{upper}$ and V$_\textrm{lower}$. 
A key achievement of \cite{TeChKiOcSmZa:00} was to explain in the language of matched asymptotic analysis how the fields associated with each of the two tangency points combine in the far field, as the associated penumbras merge. (We note also the related results obtained in \cite{MeTa:85} using the more technical machinery of microlocal analysis.) 
Without an interpretation of the divergent integral, which is not provided in \cite{TeChKiOcSmZa:00}, or in any of the related literature \cite[eqns\ (8.188)-(8.194)]{SmithThesis}, \cite[eqn\ (5.144)]{CoatsThesis}, or \cite[eqns\ (3.59) and (4.166)]{FozardThesis}, as they stand the results of \cite{TeChKiOcSmZa:00} are impossible to use for quantitative calculations; in this paper we remedy this. (Correct expressions for the transition region behaviour have been published elsewhere (see e.g.\ \cite[\S6.8]{BaKi:79} or \cite[\S13.7]{BaBu:91}), but to the best of the author's knowledge not in the cartesian coordinate system that was necessary for the global analysis in \cite{TeChKiOcSmZa:00}; our formula \rf{eqn:gDef} therefore completes the analysis of \cite{TeChKiOcSmZa:00}.)

An outline of the paper is as follows. 
In \S\ref{sec:ClassicalSoln} we state the ray tangency problem, review the classical solution in the inner Fock region, and present our new solution representation. For completeness we also detail some elementary properties of the Pekeris caret function, including some of its integral representations and its large-argument asymptotics. 
In \S\ref{sec:Matching} we use these properties to verify that the new solution representation matches correctly with the field in the outer regions, using the method of steepest descent. 
Along the way we provide the correct interpretation (in terms of the Pekeris caret function) of the divergent integral appearing in \cite[Eqn (2.63)]{TeChKiOcSmZa:00}. 
In \S\ref{sec:Conclusions} we offer some conclusions, and place the current study in the context of other ongoing work on canonical diffraction problems.
\ifpre In Appendix \ref{sec:GeneralBCs} we show how the results for the Dirichlet case (to which we restrict our attention in the main body of the paper) can be modified to deal with Neumann and Robin (mixed) boundary conditions, and in Appendix \ref{app:Airy} we define the Airy function notation used throughout the paper. 
\fi 
Finally, Appendix \ref{app:JRO} (written by J.~Ockendon) provides an alternative regularisation of the divergent integral appearing in \cite{TeChKiOcSmZa:00}, showing how the correct transition behaviour can be formally obtained from the classical Fourier integral solution representation by the method of stationary phase, a calculation that does not seem to have been published previously.

\ifpre
\else
For brevity we consider here only the case of Dirichlet boundary conditions. The corresponding formulae for Neumann and Robin (impedance) boundary conditions can be found in \cite[Appendix A]{FLP_preprint}.
\fi

\section{Problem statement and the field in the Fock region}
\label{sec:ClassicalSoln}
We consider the two-dimensional time-harmonic scattering (with time dependence $\re^{-\ri \omega t}$, $\omega>0$, assumed throughout) of an incident plane wave $\phi^i=\re^{\ri k x}$ by a smooth convex scatterer $D$.  
We seek a scattered field $\phi^s$ which solves the dimensionless Helmholtz equation (here $k>0$ is the wavenumber)
\begin{align}
\label{eqn:HE}
\pdtwo{\phi}{x}+\pdtwo{\phi}{y} + k^2 \phi =0,
\end{align}
and which, when added to the incident field, gives a total field $\phi:=\phi^i + \phi^s$ satisfying 
the Dirichlet (sound soft) boundary condition 
\begin{align}
\label{eqn:DirichletBC}
\phi=0, \qquad \textrm{on }\partial D, 
\end{align}
where $\partial D$ denotes the boundary of $D$. 

The convexity of $D$ means there are two points on $\partial D$ where the incident rays, which point in the positive $x$ direction, intersect $\partial D$ tangentially. Our focus is on the high frequency ($k\to \infty$) behaviour of the wave field in a fixed ($k$-independent) neighbourhood of one of these two tangency points.\footnote{The full wave field also involves a component associated with the second tangency point, but we do not consider this here; as was mentioned in \S\ref{sec:intro}, the far-field merging of the fields associated with the two tangency points has been described in \cite{MeTa:85,TeChKiOcSmZa:00}.} 
We denote the tangency point in question by $C$ and without loss of generality assume that it lies at the origin $(x,y)=(0,0)$. We seek a total field in the vicinity of $C$ of the form
\begin{align}
\label{eqn:HETotalField}
\phi = A\re^{\ri kx},
\end{align}
where the amplitude $A=A^i+A^s$ is a function of $x$ and $y$, with $A^i\equiv 1$ denoting the contribution from the incident field and $A^s$ denoting the scattered amplitude. (Our choice of cartesian rather than curvilinear coordinates is made to allow easier comparison with \cite{TeChKiOcSmZa:00}.)

We assume that $\partial D$ is locally parabolic near $C$, with the local curvature $\kappa$ of $\partial D$ at $C$ satisfying $0<\kappa \ll k$. In fact for ease of presentation we shall assume throughout that $\kappa=1/2$, 
so that the local form of $\partial D$ near $C$ is
\begin{align}
\label{eqn:Parabola}
y+\frac{x^2}{4} = 0.
\end{align}
Results for the general case can be obtained by replacing $x$ and $y$ respectively by $(2\kappa)^{2/3}x$ and $(2\kappa)^{1/3}y$ in all formulas from \rf{eqn:Parabola} onwards. (The case of a sharp tip where $\kappa \sim \ord{k}$ requires a different analysis; see e.g.\ \cite{EnKiTe:98,FozardThesis}).

\subsection{Classical solution representation in the Fock region}
\label{sec:FLClassicalSoln}
Our main focus is on the field in the Fock region (region I in \F\ref{fig:OuterRegions}), a neighbourhood of the point $C$ of width $\ord{k^{-1/3}}$ and height $\ord{k^{-2/3}}$. 
Following the Fock-Leontovich approach, scaling $x=k^{-1/3}\hx$, $y=k^{-2/3}\hy$, with $\hx$ and $\hy$ both $\ord{1}$, one finds that, to leading order as $k\to\infty$, the amplitude $A$ of \rf{eqn:HETotalField} (and hence also $A^s$) is a solution of the PWE
\begin{align}
\label{eqn:PWE}
2\ri \pdone{A}{\hx}+\pdtwo{A}{\hy} = 0,
\end{align}
supplemented with the boundary condition
\begin{align}
\label{eqn:FLDirichletBC}
A=0, \qquad \textrm{on } \hy+\frac{\hx^2}{4}=0.
\end{align}
As $\sqrt{\hx^2+\hy^2}\to \infty$ this inner solution must also match correctly to the inner limits of the fields in the outer regions II, III, IV, V and VI (cf.\ \S\ref{sec:Matching}).

Fourier transform methods lead to the classical solution representation (sometimes called the ``Fock formula'') 
for the scattered amplitude 
\cite[Eqn (2.15)]{TeChKiOcSmZa:00}
\begin{align}
A^s &= -\re^{-\ri(\hx \hy/2+ \hx^3/12)}\int_{-\infty}^\infty\re^{\ri \hx\sigma/2} \frac{\rA_0(\sigma)}{\rA_1(\sigma)}\rA_1(\sigma-\hn) \,\rd\sigma,
\label{eqn:FLScatClassical}
\end{align}
where $\rA_0(z):=\Ai(z)$, $\rA_1(z):=\re^{\ri 2\pi/3}\Ai(\re^{\ri 2\pi/3}z)$ are Airy 
\ifpre
functions (as defined in Appendix \ref{app:Airy}), 
\else
functions,
\fi
\begin{align}
\label{eqn:nDef}
\hn:=\hy+\hx^2/4
\end{align}
is (to leading order) the (scaled) normal distance from the observation point to the boundary\footnote{We remark that \rf{eqn:FLScatClassical} agrees with the formula presented in \cite[equation (2.80)]{OckTew:12}, provided that $2$ is corrected to $2^{1/3}$ in the denominator in \cite[equation (2.81)]{OckTew:12}. For completeness we also correct two further typographical errors in \cite{OckTew:12}: in the line before equation (2.82), and in equation (2.83) in \cite{OckTew:12}, $2\pi/3$ should be replaced by $\ri2\pi/3$.}. 
\ifpre
The Fourier transform relation \rf{eqn:AiryRelation3} then allows us to write the total amplitude $A$ as (cf.\ \cite[\S6.8]{BaKi:79})
\else
The total amplitude $A$ can be written as (cf.\ \cite[\S6.8]{BaKi:79})
\fi
\begin{align}
A &= \re^{-\ri(\hx \hy/2+ \hx^3/12)}\int_{-\infty}^\infty\re^{\ri \hx\sigma/2}\left[\rA_0(\sigma-\hn)-\frac{\rA_0(\sigma)}{\rA_1(\sigma)}\rA_1(\sigma-\hn)\right] \,\rd\sigma.
\label{eqn:FLTotClassical}
\end{align}

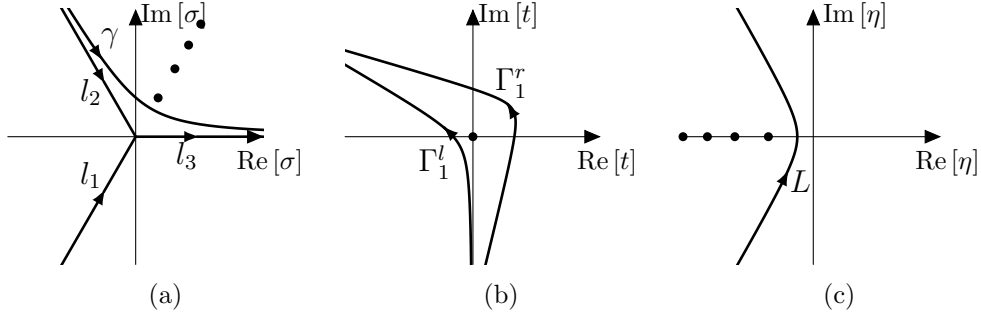
\begin{figure}[t]
\colorlet{lightgray}{black!15}
\begin{center}
\subfigure[(a)]{
\begin{tikzpicture}[line cap=round,line join=round,>=triangle 45,x=1.0cm,y=1.0cm, scale=1.7]
\def\xm{1};
\def\ym{1};
\pgfmathsetmacro{\rad}{2*sqrt(\xm^2+\ym^2)}
\def\theta{60};
\pgfmathsetmacro{\costheta}{cos(\theta)};
\pgfmathsetmacro{\sintheta}{sin(\theta)};
\pgfmathsetmacro{\costwotheta}{cos(2*\theta)};
\pgfmathsetmacro{\sintwotheta}{sin(2*\theta)};
\pgfmathsetmacro{\costhreetheta}{cos(3*\theta)};
\pgfmathsetmacro{\sinthreetheta}{sin(3*\theta)};
\pgfmathsetmacro{\cosfourtheta}{cos(4*\theta)};
\pgfmathsetmacro{\sinfourtheta}{sin(4*\theta)};
\pgfmathsetmacro{\cosfivetheta}{cos(5*\theta)};
\pgfmathsetmacro{\sinfivetheta}{sin(5*\theta)};
\pgfmathsetmacro{\cossixtheta}{cos(6*\theta)};
\pgfmathsetmacro{\sinsixtheta}{sin(6*\theta)};
\def\sc{0.15}
\def\nn{0.4}
\fill (\sc*2.34*\costheta,\sc*2.34*\sintheta) circle (1pt);
\fill (\sc*4.09*\costheta,\sc*4.09*\sintheta) circle (1pt);
\fill (\sc*5.52*\costheta,\sc*5.52*\sintheta) circle (1pt);
\fill (\sc*6.79*\costheta,\sc*6.79*\sintheta) circle (1pt);
\draw (\xm+0.05,-0.2) node {\footnotesize{$\real{\sigma}$}};
\draw [->] (-\xm,0) -- (\xm,0);
\draw [->] (0,-\ym) -- (0,\ym);
\draw (0.3,\ym-0.05) node {\footnotesize{$\im{\sigma}$}};
\draw (-0.2,0.77) node {$\gamma$};
\draw (-0.35,0.35) node {$l_2$};
\draw (-0.35,-0.3) node {$l_1$};
\draw (0.4,-0.15) node {$l_3$};
\clip (-\xm,-\ym) rectangle (\xm,\ym);
\begin{scope}[>=latex,line width=1pt]
\draw [->] (1.2*\xm*\costwotheta,1.2*\ym*\sintwotheta) -- (\xm/2.1*\costwotheta,\ym/2.1*\sintwotheta);
\draw [-] (\xm/2*\costwotheta,\ym/2*\sintwotheta) -- (0,0);
\draw [->] (1.2*\xm*\cosfourtheta,1.2*\ym*\sinfourtheta) -- (\xm/2.1*\cosfourtheta,\ym/2.1*\sinfourtheta);
\draw [-] (\xm/2*\cosfourtheta,\ym/2*\sinfourtheta) -- (0,0);
\draw [->] (0,0) -- (\xm/2,0);
\draw [-] (\xm/2.1,0) -- (\xm,0);
\end{scope}
\begin{scope}[scale=0.25,>=latex,rotate around={-30:(0,0)},line width=1]
\draw [samples=50,domain=-0.99:0.99,rotate around={90:(0,0)},xshift=0cm,yshift=0cm] plot ({1*(1+(\x)^2)/(1-(\x)^2)},{1.73*2*(\x)/(1-(\x)^2)});
\draw [<-] (-2.1,1.55) -- (-2.4,1.7);
\end{scope}
\end{tikzpicture} 
}
\subfigure[(b)]{
\begin{tikzpicture}[line cap=round,line join=round,>=triangle 45,x=1.0cm,y=1.0cm, scale=1.7]
\def\xm{1};
\def\ym{1};
\pgfmathsetmacro{\rad}{2*sqrt(\xm^2+\ym^2)}
\def\theta{60};
\pgfmathsetmacro{\costheta}{cos(\theta)};
\pgfmathsetmacro{\sintheta}{sin(\theta)};
\pgfmathsetmacro{\costwotheta}{cos(2*\theta)};
\pgfmathsetmacro{\sintwotheta}{sin(2*\theta)};
\pgfmathsetmacro{\costhreetheta}{cos(3*\theta)};
\pgfmathsetmacro{\sinthreetheta}{sin(3*\theta)};
\pgfmathsetmacro{\cosfourtheta}{cos(4*\theta)};
\pgfmathsetmacro{\sinfourtheta}{sin(4*\theta)};
\pgfmathsetmacro{\cosfivetheta}{cos(5*\theta)};
\pgfmathsetmacro{\sinfivetheta}{sin(5*\theta)};
\pgfmathsetmacro{\cossixtheta}{cos(6*\theta)};
\pgfmathsetmacro{\sinsixtheta}{sin(6*\theta)};
\begin{scope}
\clip (-\xm,-\ym) rectangle (\xm,\ym);
\begin{scope}[scale=0.125,>=latex,line width=1]
\begin{scope}[rotate around={120:(0,0)}]
\draw [samples=50,domain=-0.99:0.99,rotate around={90:(0,0)},xshift=0cm,yshift=0cm] plot ({1*(1+(\x)^2)/(1-(\x)^2)},{1.73*2*(\x)/(1-(\x)^2)});
\draw [->] (1.4,1.27) -- (1.5,1.3);
\end{scope}
\begin{scope}[rotate around={210:(0,0)},xshift=-4cm,yshift=-0.2cm,y=0.6cm]
\draw [samples=50,domain=-0.99:0.99,rotate around={0:(0,0)}] plot ({1*(1+(\x)^2)/(1-(\x)^2)},{1.73*2*(\x)/(1-(\x)^2)});
\draw [->] (1.1,-0.8) -- (1.1,-0.9);
\end{scope}
\end{scope}
\end{scope}
\draw (\xm+0.05,-0.2) node {\footnotesize{$\real{ t}$}};
\draw [->] (-\xm,0) -- (\xm,0);
\draw [->] (0,-\ym) -- (0,\ym);
\draw (0.28,\ym-0.05) node {\footnotesize{$\im{ t}$}};
\draw (0.3*\xm,0.42*\ym) node {$\Gamma_{1}^r$};
\draw (-0.3*\xm,-0.2*\ym) node {$\Gamma_{1}^l$};
\fill (0,0) circle (1pt);
\end{tikzpicture} 
}
\subfigure[(c)]{
\begin{tikzpicture}[line cap=round,line join=round,>=triangle 45,x=1.0cm,y=1.0cm, scale=1.7]
\def\xm{1};
\def\ym{1};
\pgfmathsetmacro{\rad}{2*sqrt(\xm^2+\ym^2)}
\def\theta{60};
\pgfmathsetmacro{\costheta}{cos(\theta)};
\pgfmathsetmacro{\sintheta}{sin(\theta)};
\pgfmathsetmacro{\costwotheta}{cos(2*\theta)};
\pgfmathsetmacro{\sintwotheta}{sin(2*\theta)};
\pgfmathsetmacro{\costhreetheta}{cos(3*\theta)};
\pgfmathsetmacro{\sinthreetheta}{sin(3*\theta)};
\pgfmathsetmacro{\cosfourtheta}{cos(4*\theta)};
\pgfmathsetmacro{\sinfourtheta}{sin(4*\theta)};
\pgfmathsetmacro{\cosfivetheta}{cos(5*\theta)};
\pgfmathsetmacro{\sinfivetheta}{sin(5*\theta)};
\pgfmathsetmacro{\cossixtheta}{cos(6*\theta)};
\pgfmathsetmacro{\sinsixtheta}{sin(6*\theta)};
\def\sc{0.15}
\def\nn{0.4}
\fill (\sc*2.34*\costhreetheta,\sc*2.34*\sinthreetheta) circle (1pt);
\fill (\sc*4.09*\costhreetheta,\sc*4.09*\sinthreetheta) circle (1pt);
\fill (\sc*5.52*\costhreetheta,\sc*5.52*\sinthreetheta) circle (1pt);
\fill (\sc*6.79*\costhreetheta,\sc*6.79*\sinthreetheta) circle (1pt);
\draw (\xm+0.05,-0.2) node {\footnotesize{$\real{\eta}$}};
\draw [->] (-\xm,0) -- (\xm,0);
\draw [->] (0,-\ym) -- (0,\ym);
\draw (0.33,\ym-0.05) node {\footnotesize{$\im{\eta}$}};
\draw (-0.1,-0.35) node {$L$};
\clip (-\xm,-\ym) rectangle (\xm,\ym);
\begin{scope}[scale=0.125,>=latex,rotate around={90:(0,0)},line width=1]
\draw [samples=50,domain=-0.99:0.99,rotate around={90:(0,0)},xshift=0cm,yshift=0cm] plot ({1*(1+(\x)^2)/(1-(\x)^2)},{1.73*2*(\x)/(1-(\x)^2)});
\draw [<-] (-2.1,1.55) -- (-2.4,1.7);
\end{scope}
\end{tikzpicture} 
}
\caption{Integration contours.  
}
\label{fig:gamma}
\end{center}
\end{figure}

A major drawback of the representations \rf{eqn:FLScatClassical} and \rf{eqn:FLTotClassical} is that the convergence of the integrals is very delicate. Using the well-known large argument asymptotics of the Airy \ifpre
functions (reviewed in \rf{eqn:AiryAsympt1}-\rf{eqn:AiryPrimeAsympt2}), 
\else functions, 
\fi
one finds that while both integrals are exponentially convergent as $\sigma\to+\infty$, they are defined only improperly as $\sigma\to-\infty$, because the integrand decays only algebraically (like $|\sigma|^{-1/4}$), oscillating with a phase proportional to $|\sigma|^{3/2}$. For numerical or asymptotic evaluation it is desirable to regularise the integrals. For \rf{eqn:FLTotClassical} (the expression for the total field) this can be done by deforming the path of integration onto the contour $\gamma$ illustrated in \F\ref{fig:gamma}(a); here $\gamma$ is any contour which starts at $\sigma=\re^{\ri2\pi/3}\infty$, passes below all the poles of the integrand (which lie on the line $\arg{\sigma}=\pi/3$) and ends at $\sigma=+\infty$. That the integrand decays exponentially at $\sigma=\re^{\ri2\pi/3}\infty$ can be checked by using the connection formula 
\ifpre
\rf{eqn:Connection} 
\else
\fi
\begin{align}
\label{eqn:ConnectionSpecial}
\rA_0(z)=-\rA_1(z)-\rA_2(z)
\end{align}
\ifpre
(here $\rA_2(z):=\re^{-\ri 2\pi/3}\Ai(\re^{-\ri 2\pi/3}z)$, cf.\ Appendix \ref{app:Airy}) 
\else
where $\rA_2(z):=\re^{-\ri 2\pi/3}\Ai(\re^{-\ri 2\pi/3}z)$, 
\fi
to write
\begin{align*}
\label{}
\rA_0(\sigma-\hn)-\frac{\rA_0(\sigma)}{\rA_1(\sigma)}\rA_1(\sigma-\hn)
&= \left(\frac{\rA_0(\sigma-\hn)}{\rA_1(\sigma-\hn)}-\frac{\rA_0(\sigma)}{\rA_1(\sigma)}\right)\rA_1(\sigma-\hn)\notag \\
&= \left(\frac{\rA_2(\sigma)}{\rA_1(\sigma)}-\frac{\rA_2(\sigma-\hn)}{\rA_1(\sigma-\hn)}\right)\rA_1(\sigma-\hn);
\end{align*}
we then 
\ifpre 
apply \rf{eqn:AiryAsympt1} and 
\fi
expand for large $|\sigma|$ 
to see that
\begin{align}
\label{eqn:FDifferenceTimesAiryAsymptotics}
\rA_0(\sigma-\hn)-\frac{\rA_0(\sigma)}{\rA_1(\sigma)}\rA_1(\sigma-\hn)
&= \ord{
 \dfrac{\re^{-(2/3) (\re^{-\ri 2\pi/3}\sigma)^{3/2}+\re^{\ri \pi/3}\hn(\re^{-\ri 2\pi/3}\sigma)^{1/2}}}{|\sigma|^{1/4}}},
 \notag  \\
& \vspace{-5mm}  |\sigma| \to \infty, \,\, n=\ord{1},\,\, \arg{\sigma}\in(\pi/3,\pi).
\end{align}

Regularising the expression for the scattered field \rf{eqn:FLScatClassical} is more difficult, because the integrand grows exponentially at infinity both above and below the negative real axis. The approach adopted in \cite{BaKi:79,BaBu:91} involves splitting the integral in \rf{eqn:FLScatClassical} into two integrals over $(-\infty,0)$ and $(0,\infty)$, and using the connection formula \rf{eqn:ConnectionSpecial} to write the (slowly convergent) integral over $(-\infty,0)$ as a sum of two separate integrals over $(-\infty,0)$ corresponding to the two terms arising from \rf{eqn:ConnectionSpecial}. One then deforms the integration contour \emph{differently for each of these two terms}, to arrive at the so-called ``forked contour'' representation
\begin{align}
\label{}
A^s &= \re^{-\ri(\hx \hy/2+ \hx^3/12)}
\left(\int_{l_1} \re^{\ri \hx\sigma/2}\rA_1(\sigma-\hn)\,\rd \sigma
+ \int_{l_2} \re^{\ri \hx\sigma/2}\frac{\rA_2(\sigma)}{\rA_1(\sigma)}\rA_1(\sigma-\hn) \,\rd \sigma \right. \notag \\
&\hs{50} \left. - \int_{l_3} \re^{\ri \hx\sigma/2}\frac{\rA_0(\sigma)}{\rA_1(\sigma)}\rA_1(\sigma-\hn) \,\rd \sigma\right),
\label{eqn:FLScatHorny}
\end{align}
where the contours $l_1$, $l_2$ and $l_3$ go from $\re^{-\ri2\pi/3}\infty$ to $0$, from $\re^{\ri2\pi/3}\infty$ to $0$ and from $0$ to $\infty$ respectively (see \F\ref{fig:gamma}(a)). 

\subsection{New solution representation}
\label{sec:NewSoln}
The main purpose of this paper is to report a new representation for $A^s$, equivalent to \rf{eqn:FLScatClassical} and \rf{eqn:FLScatHorny}, comprising a single rapidly convergent contour integral. 
To give some context to our new representation, we note that \rf{eqn:FLScatClassical} and \rf{eqn:FLScatHorny} both express $A^s$ as a superposition of PWE solutions of the form
\begin{align*}
\label{}
\re^{-\ri(\hx \hy/2+ \hx^3/12)}\re^{\ri \hx\sigma/2}\rA_1(\sigma-\hy-\hx^2/4), \qquad \sigma\in \C.
\end{align*}
Our new representation expresses $A^s$ as a superposition of much simpler PWE solutions
\begin{align*}
\label{}
\re^{-\ri(\hy t+ \hx t^2/2)}, \qquad t\in \C.
\end{align*}

Specifically, we claim that 
\begin{align}
\label{eqn:FLScatNew}
A^s = \int_{\Gamma_1^l} \hp(t) \re^{\ri(-\hy t-\hx t^2/2 + t^3/3)}\,\rd t,
\end{align}
where $\Gamma_1^l$ is any contour going from $-\ri\infty$ to $\re^{\ri 5\pi/6}\infty$, passing to the left of the origin (see \F\ref{fig:gamma}(b)), and $\hp(t)$ is the ``Pekeris caret function'', 
a meromorphic function of $t$ with a single simple pole at $t=0$.   
For $t\neq 0$, $\hp(t)$ can be defined by
(cf.\ the papers by Logan \cite[\S7]{Lo:59} and James \cite[pp.~25-26]{James})\footnote{In order to make formula \rf{eqn:FLScatNew} as simple as possible, we have used a slightly different normalisation to that used by Logan \cite{Lo:59} and James \cite{James}; specifically, $\hp(t)=(\ri/\sqrt{\pi})\hp_{\rm Logan}(t)=(\ri/\sqrt{\pi})\overline{\hp_{\rm James}(\overline{t})}$ and $p(t)=(\ri/\sqrt{\pi})p_{\rm Logan}(t)=(\ri/\sqrt{\pi})\overline{p_{\rm James}(\overline{t})}$. The complex conjugation relating $\hp$ to $\hp_{\rm James}$ and $p$ to $p_{\rm James}$ is due to the fact that James assumes $\re^{\ri \omega t}$ time dependence in \cite{James}.}
\begin{align}
\hp(t) = \frac{1}{2\pi\ri t} + p(t),
\qquad t\neq 0,
\label{eqn:PekerisRep1b}
\end{align}
where $p(t)$ is the ``Pekeris function'' \cite[p.~26]{James}, an entire function defined by
\begin{align}
\label{eqn:PekerisDefn}
p(t)=\frac{1}{2\pi}\left(\int_{l_2} \re^{\ri t\sigma}\frac{\rA_2(\sigma)}{\rA_1(\sigma)} \,\rd \sigma
- \int_{l_3} \re^{\ri t\sigma}\frac{\rA_0(\sigma)}{\rA_1(\sigma)} \,\rd \sigma\right), \qquad t\in\C.
\end{align}
However, as we shall review in \S\ref{sec:OtherReps}, a number of other representations for $\hp(t)$ are available. Practical computation of $\hp(t)$ is considered in \cite{Pe:87,Sa:09}.

As alluded to in \S\ref{sec:intro}, the function $\hp(t)$ is already well known to play a key role in describing the field in the penumbra region (see \S\ref{subsec:penumbra}). What formula \rf{eqn:FLScatNew} reveals, apparently for the first time, is that $\hp(t)$ is in fact also intrinsic to the description of the field in the inner Fock region. 
To expand on the historical remarks made in \S\ref{sec:intro}, the association of $\hp(t)$ and $p(t)$ with the name of Pekeris stems from the appearance of the function $p(t)$ Pekeris' 1947 paper \cite{Pe:47}\footnote{Explicitly, $p(t)=(1/12)^{1/3}\re^{\ri\pi/3}/(2\pi)\overline{F((3/2)^{2/3}\overline{t})}$, where $F(t)$ is the function defined by Pekeris in \cite[eqn (68)]{Pe:47}.}. 
But we note that $p(t)$ also appears throughout the Russian literature, for example in a 1948 paper by Fock, reproduced in \cite[\S4]{Fo:65}\footnote{Explicitly, $p(t)=\re^{\ri\pi/4}/\sqrt{\pi}g_{\rm Fock}(t)$, where $g_{\rm Fock}(t)$ is the function appearing in \cite[Eqn (4.03)]{Fo:65}.}. Hence the function $\hp(t)$ is referred to by some as a ``Fock-type integral'' \cite{Pe:87} or ``Fock scattering function'' \cite{Sa:09}.

Before reviewing some of the key properties of the function $\hp(t)$ in \S\ref{sec:OtherReps}, and providing a rigorous derivation of \rf{eqn:FLScatNew} in \S\ref{sec:Derivation}, we highlight two attractive features of our new solution representation. 
First, the integral in \rf{eqn:FLScatNew} converges rapidly: it follows from \rf{eqn:pPekerisAsympt2} below that, for fixed $(\hx,\hy)$, as $|t|\to\infty$ in either direction along the contour $L$ the integrand decays faster than $\re^{-c|t|^3}$ for every $c<1/4$. 
Second, as is clear from \rf{eqn:PekerisRep1b}, the residue of $\hp(t)$ at the pole $t=0$ is $1/(2\pi\ri)$. Hence an expression for the total amplitude $A$, equivalent to the classical representation \rf{eqn:FLTotClassical}, can be obtained from \rf{eqn:FLScatNew} by simply deforming the integration contour across the origin, giving
\begin{align}
\label{eqn:FLTotNew}
A = \int_{\Gamma_1^r} \hp(t) \re^{\ri(-\hy t-\hx t^2/2 + t^3/3)}\,\rd t,
\end{align}
where $\Gamma_1^r$ is any contour going from $-\ri\infty$ to $\re^{\ri 5\pi/6}\infty$, passing to the \emph{right} of the origin (see \F\ref{fig:gamma}(b)).

\subsection{Elementary properties of $\hp(t)$}
\label{sec:OtherReps}
In this section we document some elementary properties of the function $\hp(t)$. 
Some of these properties can be found in \cite{Lo:59,James}, but since these references do not seem to be widely known, and we use a different normalisation compared to \cite{Lo:59,James}, we restate them here for ease of reference. 

We first consider some alternative representations for $\hp(t)$, valid in different sectors of the complex plane. 
The definitions \rf{eqn:PekerisRep1b}-\rf{eqn:PekerisDefn} immediately imply the following ``forked contour'' representation, valid for $\arg{t}\in(2\pi/3,5\pi/3)$, which will be used in the derivation of \rf{eqn:FLScatNew} in \S\ref{sec:Derivation}:
\begin{align}
\hp(t) &= \frac{1}{2\pi}\left(\int_{l_1} \re^{\ri t\sigma}\,\rd \sigma
+ \int_{l_2} \re^{\ri t\sigma}\frac{\rA_2(\sigma)}{\rA_1(\sigma)} \,\rd \sigma
- \int_{l_3} \re^{\ri t\sigma}\frac{\rA_0(\sigma)}{\rA_1(\sigma)} \,\rd \sigma\right),
\notag\\& \hs{75} \arg{t}\in(2\pi/3,5\pi/3).
\label{eqn:PekerisRep1c}
\end{align}
Deforming the contours for the first two terms in \rf{eqn:PekerisRep1c} onto the negative real axis and applying \rf{eqn:ConnectionSpecial} gives a Fourier-type representation valid for $\im{t}<0$:
\begin{align}
\label{eqn:PekerisRep2}
\hp(t) = -\frac{1}{2\pi}\int_{-\infty}^\infty \re^{\ri t\sigma}\frac{\rA_0(\sigma)}{\rA_1(\sigma)} \,\rd \sigma,
\quad \im{t}<0.
\end{align}
By a further contour deformation one can obtain a representation valid for $ \arg{t}\in(-2\pi/3,\pi/3)$ (here the convergence at $\sigma=+\infty$ is ensured because of the exponential decay of $\rA_0(\sigma)/\rA_1(\sigma)$; at $\sigma=\re^{\ri 2\pi/3}\infty$, $\rA_0(\sigma)/\rA_1(\sigma)$ is $\ord{1}$, and the convergence relies on the exponential decay of the factor $\re^{\ri t\sigma}$):
\begin{align}
\label{eqn:PekerisRep1}
\hp(t) = -\frac{1}{2\pi}\int_{\gamma} \re^{\ri t\sigma}\frac{\rA_0(\sigma)}{\rA_1(\sigma)} \,\rd \sigma.
\quad \arg{t}\in(-2\pi/3,\pi/3).
\end{align}
Finally, integrating by parts using the Wronskian relation 
\ifpre 
\rf{eqn:AiryWronskian}, 
\else
\cite[\S9.2]{DLMF}\footnote{
We warn the reader that the Wronskian formulas in \cite[p.~405]{BaBu:91} are incorrect: instead of $\{2,\ri,\ri,1\}$ they should read $\{2\ri,-1,-1,-1\}$.}, 
\fi
and changing variable to $\eta=\re^{\ri 2\pi/3}\sigma$,
gives a representation valid for all $t\in \C\setminus\{0\}$:
\begin{align}
\hp(t) = -\frac{1}{4\pi^{2}t}\int_{L} \frac{\re^{\re^{-\ri\pi/6}t\eta}}{\rA_0(\eta)^2} \,\rd \eta,\qquad t\in\C\setminus\{0\}.
\label{eqn:pPekeris}
\end{align}
Here $L$ is any contour going from $\re^{-\ri 2\pi/3}\infty$ to $\re^{\ri 2\pi/3}\infty$, passing to the right of all the poles of the integrand, at the zeros $\eta_n$, $n=0,1,2,\ldots$, of the Airy function, all of which lie on the negative real axis (see \F\ref{fig:gamma}(c)). 


Deforming the contour in \rf{eqn:pPekeris} to wrap around the poles of the integrand on the negative real axis gives the residue series representation 
\begin{align}
\label{eqn:PekerisResidue}
\hp(t) = \frac{\re^{-\ri 2\pi/3}}{2\pi} \sum_{n=0}^\infty \frac{\re^{\re^{-\ri\pi/6}t\eta_{n}}}{\rA_0'(\eta_{n})^2}, \qquad \arg{t}\in(-\pi/3,2\pi/3).
\end{align}
The large argument asymptotics 
are governed by the first term, i.e.
\begin{align}
\label{eqn:pPekerisAsympt1}
\hp(t) \sim \frac{\re^{-\ri 2\pi/3}}{2\pi} \frac{\re^{\re^{-\ri\pi/6}t\eta_{0}}}{\rA_0'(\eta_{0})^2}, \qquad |t|\to\infty,\,\arg{t}\in(-\pi/3,2\pi/3),
\end{align}
the remainder being exponentially small compared to this term. 

The large argument asymptotics of $\hp(t)$ in the sector $(2\pi/3,5\pi/3)$ can be obtained from \rf{eqn:pPekeris} using the method of steepest descent. Using the large argument asymptotics of the Airy 
\ifpre
function (cf.~\rf{eqn:AiryAsympt1} and \rf{eqn:AiryPrimeAsympt1}),
\else 
function,
\fi
away from the negative real axis we can approximate the integrand in \rf{eqn:pPekeris} for large $|\eta|$ by
\begin{align*}
\label{}
4\pi \eta^{1/2} \re^{\re^{-\ri\pi/6}t\eta+(4/3)\eta^{3/2}}\left(1+ \ord{\frac{1}{|\eta|^{3/2}}} \right).
\end{align*}
The terms in the exponent are comparable in magnitude when $\eta\sim \ord{|t|^2}$; rescaling $\eta=|t|^2\zeta$ gives
\begin{align*}
\label{}
\hp(t) \sim -\frac{|t|^2}{\pi \vartheta}\int_L  \zeta^{1/2} \re^{|t|^3(\re^{-\ri\pi/6}\vartheta\zeta+(4/3)\zeta^{3/2})} \left(1+ \ord{\frac{1}{|t|^3|\zeta|^{3/2}}} \right) \,\rd \zeta,
\end{align*}
where $t=|t|\vartheta$, $|\vartheta|=1$. The phase is stationary when
\begin{align}
\label{eqn:StatPoint}
\zeta^{1/2}=-\re^{-\ri\pi/6}\vartheta/2,
\end{align}
and, recalling that the square root takes its principal branch, a solution of \rf{eqn:StatPoint} exists for all $\arg{t}=\arg{\vartheta}\in(2\pi/3,5\pi/3)$, namely 
\begin{align*}
\label{}
\zeta_* = \re^{-\ri\pi/3}\vartheta^2/4.
\end{align*}
Deforming $\gamma$ onto the steepest descent contour through $\zeta_*$ then 
gives
\begin{align}
\label{eqn:pPekerisAsympt2}
\hp(t)\sim \frac{\sqrt{-t}}{2\sqrt{\pi}}\re^{-\ri (t^3/12-\pi/4)} \left(1+ \ord{\frac{1}{|t|^3}} \right), \quad |t|\to\infty,\,\arg{t} \in(2\pi/3,5\pi/3).
\end{align}
Higher order terms (all the way up to $\ord{1/|t|^{30}}$) can be found in \cite[\S10]{Lo:59}. 

\subsection{Derivation of the new solution representation}
\label{sec:Derivation}
The key to deriving \rf{eqn:FLScatNew} is the observation that the Airy function factor $\rA_1(\sigma-\hn)$ appearing in \rf{eqn:FLScatClassical} and \rf{eqn:FLScatHorny} can be expressed in integral form using the identity 
\ifpre 
(cf.\ \rf{eqn:AiryIdentity})
\else
(cf.\ \cite[Eqn (9.5.4)]{DLMF}, for more details see \cite[Appendix A]{FLP_preprint})
\fi
\begin{align}
\label{eqn:AiryIdentitySpecial}
\re^{\ri \hx\sigma/2}\re^{-\ri(\hx\hy/2+\hx^3/12)}\rA_1(\sigma-\hn) =\frac{1}{2\pi} \int_{\Gamma_1}\re^{\ri \sigma t}\re^{\ri(-\hy t-\hx t^2/2 + t^3/3)} \, \rd t,
\end{align}
where 
\ifpre (cf.\ \F\ref{Sectors})  
\else 
\fi
$\Gamma_1$ is any contour going from $-\ri\infty$ to $\re^{\ri 5\pi/6}\infty$. 
In the light of \rf{eqn:PekerisRep2}, the new representation \rf{eqn:FLScatNew} can be obtained formally from the classical Fourier integral solution \rf{eqn:FLScatClassical} by first applying \rf{eqn:AiryIdentitySpecial}, and then switching the order of integration in the resulting double integral. However, to make this formal calculation rigorous, we start instead from the regularised solution \rf{eqn:FLScatHorny}. 
Substituting \rf{eqn:AiryIdentitySpecial} into \rf{eqn:FLScatHorny} gives
\begin{align}
A^s &= \frac{1}{2\pi}\left(\int_{l_1}  \int_{\Gamma_1}\re^{\ri \sigma t}f(t;\hx,\hy)\, \rd t\,\rd \sigma
+ \int_{l_2}  \int_{\Gamma_1}\re^{\ri \sigma t}f(t;\hx,\hy) \frac{\rA_2(\sigma)}{\rA_1(\sigma)} \, \rd t \,\rd \sigma\right.\notag\\
& \hs{45} \left.- \int_{l_3}  \int_{\Gamma_1}\re^{\ri \sigma t}f(t;\hx,\hy)\frac{\rA_0(\sigma)}{\rA_1(\sigma)} \, \rd t \,\rd \sigma\right),
\label{eqn:DoubleInt}
\end{align}
where $f(t;\hx,\hy)=\re^{\ri(-\hy t-\hx t^2/2 + t^3/3)}$. 
Provided that $\Gamma_1$ lies entirely in the half plane $\arg{t}\in(2\pi/3,5\pi/3)$ (in particular, passing to the left of the origin), we claim that the order of integration can be switched in the double integrals appearing in \rf{eqn:DoubleInt}, which produces the new representation \rf{eqn:FLScatNew} with the function $\hp(t)$ represented by \rf{eqn:PekerisRep1c}.

To justify the switch in integration order we consider a specific choice of $\Gamma_1$. 
Given $T>0$, let $\Gamma_1= \Gamma_{1,1}\cup \Gamma_{1,2} \cup \Gamma_{1,3}$, where $\Gamma_{1,1}$ goes from $-\ri\infty$ to $-\ri T$ along the imaginary $t$-axis, $\Gamma_{1,2}$ goes clockwise around the circle $|t|=T$ from $-\ri T$ to $\re^{\ri 5\pi/6}T$, and $\Gamma_{1,3}$ goes from $\re^{\ri 5\pi/6}T$ to $\re^{\ri 5\pi/6}\infty$ along the ray $\arg t= 5\pi/6$. 
The right-hand side of \rf{eqn:DoubleInt} can then be written as a sum of $9$ double integrals. Switching integration order is permitted by Fubini's theorem provided each of these integrals is absolutely convergent. To show this, fix $(\hx,\hy)$ and choose $T=T(\hx,\hy)>0$ large enough that $|f(t;\hx,\hy)|\leq \re^{-w^3/6}$ on $\Gamma_{1,1}$ and $\Gamma_{1,3}$, where $w=|t|$. Then, noting that there exists a constant $C$ such that $|\rA_2(\sigma)/\rA_1(\sigma)|\leq C\re^{-(4/3)s^{3/2}}$ on $l_2$, where $s=|\sigma|$, and similarly for $|\rA_0(\sigma)/\rA_1(\sigma)|$  on $l_3$, the absolute convergence of the integrals involving $\Gamma_{1,1}$ and $\Gamma_{1,3}$ follows from the finiteness of the integrals
\begin{align*}
\label{}
\int_0^\infty \int_T^\infty \re^{-sw/2-w^3/6} \, \rd w\, \rd s,\\
\int_0^\infty \int_T^\infty \re^{-sw/2 -(4/3)s^{3/2}-w^3/6} \, \rd w\, \rd s,\\
\int_0^\infty \int_T^\infty \re^{sw -(4/3)s^{3/2}-w^3/6} \, \rd w\, \rd s.
\end{align*}
The finiteness of the first two is obvious, for the third we note that 
\begin{align*}
\label{}
\int_0^\infty \re^{sw-(4/3)s^{3/2}} \, \rd s \sim \sqrt{\pi w}\re^{w^3/9}, \qquad w\to\infty.
\end{align*}
For the integrals involving $\Gamma_{1,2}$ we note that $|f(t;\hx,\hy)|\leq C=C(T)=C(\hx,\hy)$ on $\Gamma_{1,2}$, so after trivial estimation of the $t$ integral, absolute convergence follows from the finiteness of the following one-dimensional integrals:
\begin{align*}
\label{}
\int_0^\infty \re^{-sT/2-(4/3)s^{3/2}} \, \rd s, \qquad
\int_0^\infty \re^{sT-(4/3)s^{3/2}} \, \rd s.
\end{align*}


\section{Matching to the outer regions}
\label{sec:Matching}
We now apply the steepest descent method, combined with the large argument approximations of $\hp(t)$ presented in \S\ref{sec:OtherReps}, to show how our new solution representations \rf{eqn:FLScatNew} and \rf{eqn:FLTotNew} can be systematically matched to the field in the outer regions in \F\ref{fig:OuterRegions}. An overview of the general picture is as follows: when matching out to the illuminated region (i.e. $x<0$, or $x>0$ with $y/x\gg k^{-1/3}$) the main contribution to the integral comes from a saddle point on the negative real $t$-axis, where the approximation \rf{eqn:pPekerisAsympt2} holds and $\hp(t)$ is oscillating. As we approach the penumbra ($x>0$ with $y/x= \ord{k^{-1/3}}$) this saddle point approaches the pole in $\hp(t)$ at $t=0$, and the interaction between the saddle point and the pole is what ``switches off'' the incident wave across the penumbra. In the creeping wave region in the deep shadow we see a pair of saddle points close to the positive real $t$-axis, where the approximation \rf{eqn:pPekerisAsympt1} holds and $\hp(t)$ is exponentially small. These saddle points coalesce on the boundary curve itself, allowing recovery of the familiar Airy function description of the creeping field. 

\subsection{Matching to the illuminated region}
We first consider the matching from the inner Fock region I out to the illuminated region VI. In this case the main contribution to \rf{eqn:FLScatNew} comes from a saddle point on the negative real $t$-axis, at an $\ord{k^{1/3}}$ distance from the origin. We write \rf{eqn:FLScatNew} in the outer variables $x=k^{-1/3}\hx$, $y=k^{-2/3}\hy$, and consider the behaviour of the resulting expression as $k\to\infty$. Rescaling $t=k^{1/3}\tau$ and applying the large argument approximation \rf{eqn:pPekerisAsympt2} gives
\begin{align*}
\label{}
A^s \sim \frac{k^{1/2}\re^{\ri \pi/4}}{2\sqrt{\pi}} \int_{\Gamma_1^l}(-\tau)^{1/2} \re^{\ri k(-y\tau-x\tau^2/2+\tau^3/4)}\left(1+\ord{\frac{1}{k|\tau|^3}}\right) \,\rd \tau, \qquad k\to\infty.
\end{align*}
There are two real saddle points located at
\begin{align*}
\label{}
\tau_\pm = (2/3)(x\pm\sqrt{x^2+3y}),
\end{align*}
and we note that in the propagation domain ($y>-x^2/4$) these are always distinct. Since we are integrating along $\Gamma_1^l$, when we deform to the steepest descent path we pass only through the left-most saddle, $\tau_-$, giving
\begin{align}
\label{eqn:AsMatchingIlluminated}
A^s\sim \frac{1}{\sqrt{3}}\left(1-\frac{x}{\sqrt{x^2+3y}} \right)^{1/2}\re^{\ri k (4/27)\left(-x^3-(9/2)xy+(x^2+3y)^{3/2}\right)},
\end{align}
which can be shown to match correctly with the inner limit of the specularly reflected wave (cf.\ e.g.\ \cite[\S2.1]{TeChKiOcSmZa:00}) as one moves from region VI into region I.

\subsection{Matching to the penumbra and the interpretation of \cite[eqn (2.63)]{TeChKiOcSmZa:00}}
\label{subsec:penumbra}
As the observation point $(x,y)$ approaches the geometrical shadow boundary, the saddle point $\tau_-$ approaches the origin. Precisely, when $x>0$, $y>0$ and $y/x^2\ll 1$, we have $\tau_-\sim -y/x + \ord{y^2/x^3}$. So in particular once $y/x=\ord{k^{-1/3}}$ (i.e.,\ we are in the penumbra), the above analysis becomes invalid, because the approximation \rf{eqn:pPekerisAsympt2} for $\hp(t)$ no longer applies near the saddle point $\tau=\tau_-$. In this case we introduce the penumbra variable $\ty=k^{1/3}y$, go back to \rf{eqn:FLScatNew} and write it in terms of the variables $x=k^{-1/3}\hx$, $\ty=k^{-1/3}\hy$, leaving the integration variable $t$ unscaled, to obtain
\begin{align}
\label{eqn:FLScatNewTransition}
A^s = \int_{\Gamma_1^l}\re^{\ri t^3/3} \hp(t)  \re^{\ri k^{1/3}(-\ty t-x t^2/2)}\,\rd t,
\end{align}
which has a single saddle point at 
\begin{align}
\label{eqn:PenumbraSaddle}
t = -\frac{\ty}{x}.
\end{align} 
Deform the integration contour onto the steepest descent contour passing through the saddle point. If $\ty<0$ this requires us to cross the pole at $t=0$, so that in \rf{eqn:FLScatNewTransition} ${\Gamma_1^l}$ is replaced by ${\Gamma_1^r}$ and $A^s$ picks up a residue contribution of $-1$ (equivalently, in this case we can view our integral as representing the total field $A$, cf.\ \rf{eqn:FLTotNew}). Expanding the phase in \rf{eqn:FLScatNewTransition} around \rf{eqn:PenumbraSaddle} as 
\[-\ty t-x t^2/2=-(x/2)(t+\ty/x)^2+\ty^2/(2x),\]
the main contribution to \rf{eqn:FLScatNewTransition} will come from a neighbourhood of the saddle point of size $\ord{k^{-1/6}x^{-1/2}}$. Thus, if $|\ty|/\sqrt{x}\gg k^{-1/6}$ (i.e.\ we are in region V$_\textrm{upper}$ or V$_\textrm{lower}$) the standard steepest descent method gives
\begin{align}
\label{eqn:FLScatNewTransitionApprox}
A^s\sim -\Heaviside(-\ty) + \frac{\re^{\ri k^{1/3}\ty^2/(2 x)} }{k^{1/6}\sqrt{x}}g\left(\frac{\ty}{x}\right), \qquad |\ty|/\sqrt{x}\gg k^{-1/6},
\end{align}
where 
$\Heaviside(z)=0$ for $z<0$ and $\Heaviside(z)=1$ for $z>0$, and
\begin{align}
\label{eqn:gDef}
g(\xi) = \sqrt{2\pi}\re^{\ri 3\pi/4}\re^{-\ri \xi^3/3}\, \hp(-\xi).
\end{align}
We pause here to remark that formula \rf{eqn:gDef} provides the correct interpretation of the divergent integral appearing in \cite[eqn (2.63)]{TeChKiOcSmZa:00}. 
The results of \S\ref{sec:OtherReps} confirm that the expression in \cite[eqn (2.63)]{TeChKiOcSmZa:00} is \emph{formally} correct, in the sense that it agrees with \rf{eqn:gDef} if the Pekeris caret function is represented by the Fourier-type integral \rf{eqn:PekerisRep2}. But unfortunately the integral in \rf{eqn:PekerisRep2} is not convergent for real arguments; instead one should use either \rf{eqn:PekerisRep1b} or \rf{eqn:pPekeris}.


Returning to our steepest descent analysis, when $|\ty|/\sqrt{x}\ll k^{-1/6}$ (i.e.\ we are in region IV) the analysis leading to \rf{eqn:FLScatNewTransitionApprox} fails; the saddle point and the pole interact. Writing $\ty=k^{-1/6}\cy$, the leading order behaviour of \rf{eqn:FLScatNewTransition} is now
\begin{align*}
\label{}
A^s \sim -\Heaviside(-\cy) + \frac{\re^{\ri \cy^2/(2x)}}{2\pi \ri}\int_{-\infty}^\infty\frac{\re^{-\tau^2}\,\rd \tau}{\tau-\tau_*}=\Fr\left(-\dfrac{\cy}{\sqrt{2x}}\right), \qquad |\ty|/\sqrt{x}\ll k^{-1/6},
\end{align*}
where $\tau_* =\re^{-\ri 3\pi/4}\cy/\sqrt{2x}$, $\Fr(z):=(\re^{-\ri \pi/4}/\sqrt{\pi})\int_z^\infty \re^{\ri\zeta^2 \,\rd\zeta}$ is the Fresnel integral, and we have used the identity \cite[(7.5.2), (7.5.9), (7.7.2)]{DLMF}
\begin{align}
\label{eqn:FLScatNewFresnelApprox}
\int_{-\infty}^\infty\frac{\re^{-\tau^2}\,\rd \tau}{\tau-\tau_*} = 2\pi \ri \re^{-\tau_*^2}\left(\Fr(\re^{-\ri \pi/4}\tau_*) - \Heaviside(-\im{\tau_*})\right).
\end{align}
Combining \rf{eqn:FLScatNewFresnelApprox} with \rf{eqn:FLScatNewTransitionApprox} we gives a complete description of the field in the penumbra. 
With $0<x\leq \ord{1}$ we have
\begin{align*}
\label{}
A \sim 
\left\lbrace 
\begin{array}{ll}
1+ \dfrac{\re^{\ri k^{1/3}\ty^2/(2 x)} }{k^{1/6}\sqrt{x}}g\left(\dfrac{\ty}{x}\right), & \dfrac{\ty}{x}=\ord{1}, \,k^{-1/6} \ll \dfrac{\ty}{\sqrt{x}}\,\,(\textrm{Region V}_{\textrm{upper}}),\vs{2}\\
\Fr\left(-\dfrac{\cy}{\sqrt{2x}}\right), & \dfrac{\cy}{\sqrt{x}}= \ord{1}\,\,(\textrm{Region IV}),\vs{2}\\
\dfrac{\re^{\ri k^{1/3}\ty^2/(2 x)} }{k^{1/6}\sqrt{x}}g\left(\dfrac{\ty}{x}\right), & -\dfrac{\ty}{x}=\ord{1}, \,k^{-1/6} \ll -\dfrac{\ty}{\sqrt{x}}\,\,(\textrm{Region V}_{\textrm{lower}}).
\end{array}
\right.
\end{align*}
That the solutions in these three regions correctly match with each other is easily verified using the fact that $\Fr(z)\sim \re^{\ri z^2}\re^{\ri \pi/4}/(2\sqrt{\pi} z)$, $z\to\infty$, and $\Fr(z)\sim 1+ \re^{\ri z^2}\re^{\ri \pi/4}/(2\sqrt{\pi} z)$, $z\to-\infty$. 

A uniform approximation, valid across the whole penumbra (i.e. in regions IV, V$_\textrm{upper}$ and V$_\textrm{lower}$), can then be obtained by summing the approximations in the different regions and subtracting their common parts. The result is that (cf.\ the results in \cite[(eqns (6.8.2)-(6.8.6)]{BaKi:79} and \cite[(eqns (13.7.7)-(13.7.11)]{BaBu:91} and \cite[eqns (5)-(6)]{BuLy:87}, which use a different coordinate system)
\begin{align}
\label{eqn:uniformPenumbra}
A^s\sim \Fr\left(-\dfrac{\cy}{\sqrt{2x}}\right) +  \dfrac{\re^{\ri k^{1/3}\ty^2/(2 x)} }{k^{1/6}\sqrt{x}}\tg\left(\dfrac{\ty}{x}\right), \qquad  \dfrac{\ty}{x}=\ord{1},
\end{align}
where 
\begin{align}
\tg(\xi) &= \sqrt{2\pi}\re^{\ri 3\pi/4}\re^{-\ri \xi^3/3}\left(\hp(-\xi)+\frac{\re^{\ri \xi^3/3}}{2\pi\ri\xi} \right)\notag\\
&= \sqrt{2\pi}\re^{\ri 3\pi/4}\re^{-\ri \xi^3/3}\left(p(-\xi)+\frac{\re^{\ri \xi^3/3}-1}{2\pi\ri\xi} \right).
\label{eqn:tgDef}
\end{align}
Note that while $g(\xi)$ has a pole at $\xi=0$, the function $\tg(\xi)$ does not.

\subsection{Matching to the creeping wave region}
Finally, we consider the matching of our new solution representation to region II (the ``Airy layer'', or creeping wave region). We now start from the expression \rf{eqn:FLTotNew} for the total field. Writing this in terms of the variables $x=k^{-1/3}\hx$ and $\hn=\hy+\hx^2/4$ (cf.\ \rf{eqn:nDef}), scaling $t=k^{1/3}T$, and recalling \rf{eqn:pPekerisAsympt1}, we find that for $\arg{T}\in (-\pi/3,2\pi/3)$ the leading order behaviour of the integrand in \rf{eqn:FLTotNew} as $k\to\infty$ is proportional to
\begin{align*}
\label{}
\exp{\ri k\left[\left(\frac{x^2}{4}- k^{-2/3}(\hn + \re^{\ri\pi/3}\eta_{0})\right)T - \frac{x}{2}T^2 +\frac{T^3}{3} \right]},
\end{align*}
and there are saddle points at
\begin{align*}
\label{}
T_\pm = \frac{x}{2} \pm k^{-1/3} \sqrt{n+ \re^{\ri\pi/3}\eta_{0}},
\end{align*}
which coalesce as $k\to\infty$ near the point $T = x/2$. Localising the integral around this point, expanding the phase, then rescaling, one finds that to leading order \rf{eqn:FLTotNew} is proportional to 
\begin{align*}
\label{}
\re^{-\ri k (xy/2+x^3/12)}\re^{\ri k^{1/3} \re^{\ri\pi/3}\eta_{0} x}\int_{\Gamma_1} \re^{\ri[-(\hn+\re^{\ri\pi/3}\eta_{0})s+s^3/3]}\,\rd s,
\end{align*}
where $s=t-k^{1/3}x/2$. (Note that in this regime the pole in $\hp(t)$ at $t=0$ plays no role, which allows us remove the superscript $^r$ from $\Gamma_1$). Re-inserting the constant of proportionality, 
\ifpre
and recalling \rf{eqn:AjDef}, 
\fi
we find that
\begin{align}
\label{eqn:MatchingCreeping}
A\sim \frac{\re^{-\ri k (xy/2+x^3/12)}\re^{\ri k^{1/3} \re^{\ri\pi/3}\eta_{0} x}}{\rA_0'(\eta_{0})^2}\rA_0(\eta_{0}+\re^{-\ri\pi/3}\hn),
\end{align}
which correctly matches the inner limit of the creeping field (cf.\ \cite[eqn (2.36)]{TeChKiOcSmZa:00}) as one moves from region II into region I. 
We note that the shed creeping ray field in the deep shadow region III away from the boundary can then be determined by matching back into region II, exactly as in \cite[\S2.4]{TeChKiOcSmZa:00}, but we do not reproduce the details here.

\section{Conclusion}
\label{sec:Conclusions}
In this paper we have presented a new solution representation for the field in the vicinity of a tangency point between an incoming ray field and the scatterer boundary, in the case of two-dimensional scalar wave scattering by a smooth convex obstacle. Our new representation takes the form of a complex contour integral of a meromorphic function (the Pekeris caret function) mutliplied by an exponential factor with polynomial exponent. 
It can be shown that contour integrals of this type can describe other classical ``thin-layer'' wave phenomena such as creeping waves, whispering gallery waves, and the field in the vicinity of caustics (for a unified description of these thin-layer phenomena from a slightly different perspective see the review article \cite{OckTew:12}). Our hope is that the further study of such integrals 
may provide a methodology with which to attack hitherto unsolved canonical problems in diffraction theory, for example the concave-convex transition at an inflection point on a boundary (studied previously e.g.\ in \cite{Pop:79,BabSmy:86,Kaz:03}).

\bibliography{inflection_bib}
\bibliographystyle{elsarticle-num}
\appendix
\ifpre
\section{Other boundary conditions}
\label{sec:GeneralBCs}
In this appendix we indicate how the results presented above can be generalised to the case of Robin (impedance) and Neumann (sound hard) boundary conditions. The formulas now also involve the derivative of the Airy function, but since the large argument asymptotics of $\Ai$ and $\Ai'$ are essentially the same (at least in terms of the arguments of the exponential/sinusoidal factors, cf.\ \rf{eqn:AiryAsympt1}-\rf{eqn:AiryPrimeAsympt2}), all of the analysis of the previous sections carries through \emph{mutatis mutandis}. We therefore simply list the formulas obtained, for easy reference, providing commentary where required. In \S\ref{sec:Robin} we state the formulas for the Robin case; the corresponding results for the Neumann case are given in \S\ref{sec:Neumann}. Throughout this section the equation labels (R.m.n) and (N.m.n) indicate respectively the Robin (R) and Neumann (N) versions of equation (m.n) from the Dirichlet case. 


\subsection{Robin case}
\label{sec:Robin}
We consider the general boundary condition
\begin{align}
\label{eqn:GeneralBC}
\pdone{\phi}{\bn}=\mu\phi, \qquad \textrm{on }\partial D,
\tagR{eqn:DirichletBC}
\end{align}
where $\bn$ is the outward unit normal vector to $\partial D$ and $\mu$ is a constant describing the scattering properties of $\partial D$. When $\mu\in\C\setminus\{0\}$, \rf{eqn:GeneralBC} represents an impedance boundary condition (modelling an absorbing boundary if $\im{\mu}>0$). We note that when $\mu=\infty$ \rf{eqn:GeneralBC} corresponds formally to the Dirichlet (sound soft) boundary condition \rf{eqn:DirichletBC} considered previously, and when $\mu=0$ to the Neumann (sound hard) boundary condition considered in \S\ref{sec:Neumann} below. 

Scaling $\mu = k^{2/3}\hmu$ (in the general case $\kappa\neq 1/2$ we would scale $\mu = (2\kappa)^{1/3}k^{2/3}\hmu$), the leading order approximation of \rf{eqn:GeneralBC} gives the boundary condition
\begin{align}
\label{eqn:FLImpBC}
\pdone{A}{\hy}+\left(\frac{\ri\hx}{2}-\hmu\right)A=0, \qquad \textrm{on } \hy+\frac{\hx^2}{4}=0,
\tagR{eqn:FLDirichletBC}
\end{align}
in the Fock region (region I), and the classical solution representation is (cf.\ \cite[eqn (4.35)]{Fo:65})
\begin{align}
\label{eqn:FLScatClassicalImp}
\tagR{eqn:FLScatClassical}
A^s &= -\re^{-\ri(\hx \hy/2+ \hx^3/12)}\int_{-\infty}^\infty\re^{\ri \hx\sigma/2} \left(\frac{\hmu \rA_0(\sigma)-\rA_0'(\sigma)}{\hmu \rA_1(\sigma) - \rA_1'(\sigma)}\right)\rA_1(\sigma-\hn) \,\rd\sigma,\\
A &= \re^{-\ri(\hx \hy/2+ \hx^3/12)}\int_{-\infty}^\infty\re^{\ri \hx\sigma/2}\left[\rA_0(\sigma-\hn)-\left(\frac{\hmu \rA_0(\sigma)-\rA_0'(\sigma)}{\hmu \rA_1(\sigma) - \rA_1'(\sigma)}\right)\rA_1(\sigma-\hn)\right] \,\rd\sigma.
\label{eqn:FLTotClassicalImp}
\tagR{eqn:FLTotClassical}
\end{align}
The ``forked contour'' representation is
\begin{align}
\label{}
A^s &= \re^{-\ri(\hx \hy/2+ \hx^3/12)}
\left(\int_{l_1} \re^{\ri \hx\sigma/2}\rA_1(\sigma-\hn)\,\rd \sigma
+ \int_{l_2} \re^{\ri \hx\sigma/2}\frac{\hmu\rA_2(\sigma)-\rA_2'(\sigma)}{\hmu\rA_1(\sigma)-\rA_1'(\sigma)}\rA_1(\sigma-\hn) \,\rd \sigma \right. \notag \\
&\hs{50} \left. - \int_{l_3} \re^{\ri \hx\sigma/2}\frac{\hmu\rA_0(\sigma)-\rA_0'(\sigma)}{\hmu\rA_1(\sigma)-\rA_1'(\sigma)}\rA_1(\sigma-\hn) \,\rd \sigma\right)
\label{eqn:FLScatHornyImp}
\tagR{eqn:FLScatHorny}
\end{align}
and our new solution representation is
\begin{align}
\label{eqn:FLScatNewImp}
A^s = \int_{\Gamma_1^l} \hV(t,\hmu) \re^{\ri(-\hy t-\hx t^2/2 + t^3/3)}\,\rd t,
\tagR{eqn:FLScatNew}
\end{align}
\begin{align}
\hV(t,\hmu) &= \frac{1}{2\pi}\left(\int_{l_2} \re^{\ri t\sigma}\,\rd \sigma
+ \int_{l_2} \re^{\ri t\sigma}\frac{\hmu\rA_2(\sigma)-\rA_2'(\sigma)}{\hmu\rA_1(\sigma)-\rA_1'(\sigma)} \,\rd \sigma
- \int_{l_3} \re^{\ri t\sigma}\frac{\hmu\rA_0(\sigma)-\rA_0'(\sigma)}{\hmu\rA_1(\sigma)-\rA_1'(\sigma)} \,\rd \sigma\right),\notag\\
& \hs{75} \arg{t}\in(-2\pi/3,\pi/3).
\label{eqn:PekerisRep1cImp}
\tagR{eqn:PekerisRep1c}
\end{align}
\begin{align}
\hV(t,\hmu) = \frac{1}{2\pi\ri t} + V(t,\hmu),
\qquad t\neq 0,
\label{eqn:PekerisRep1bImp}
\tagR{eqn:PekerisRep1b}
\end{align}
\begin{align}
V(t,\hmu)=\frac{1}{2\pi}\left(\int_{l_2} \re^{\ri t\sigma}\frac{\hmu\rA_2(\sigma)-\rA_2'(\sigma)}{\hmu\rA_1(\sigma)-\rA_1'(\sigma)} \,\rd \sigma
- \int_{l_3} \re^{\ri t\sigma}\frac{\hmu\rA_0(\sigma)-\rA_0'(\sigma)}{\hmu\rA_1(\sigma)-\rA_1'(\sigma)} \,\rd \sigma\right)
\label{eqn:PekerisDefnImp}
\tagR{eqn:PekerisDefn}
\end{align}
\begin{align}
A = \int_{\Gamma_1^r} \hV(t,\hmu) \re^{\ri(-\hy t-\hx t^2/2 + t^3/3)}\,\rd t,
\label{eqn:FLTotNewImp}
\tagR{eqn:FLTotNew}
\end{align}
\begin{align}
\hV(t,\hmu) = -\frac{1}{2\pi}\int_{-\infty}^\infty \re^{\ri t\sigma}\frac{\hmu\rA_0(\sigma)-\rA_0'(\sigma)}{\hmu\rA_1(\sigma)-\rA_1'(\sigma)} \,\rd \sigma,
\quad \im{t}<0.
\label{eqn:PekerisRep2Imp}
\tagR{eqn:PekerisRep2}
\end{align}
\begin{align}
\hV(t,\hmu) = -\frac{1}{4\pi^{2}t}\int_{L} \frac{(\hmu^2+\re^{\ri\pi/3}\eta)\re^{\re^{-\ri\pi/6}t\eta}}{\left(\hmu\rA_0(\eta)+\re^{-\ri\pi/3}\rA_0'(\eta)\right)^2} \,\rd \eta, 
\label{eqn:pPekerisImp}
\tagR{eqn:pPekeris}
\end{align}
As before, $L$ is any contour in the complex $\eta$-plane starting at infinity with $\arg{\eta}=-2\pi/3$ and ending at infinity with $\arg{\eta}=2\pi/3$, passing to the right of all the (countably many) poles of the integrand. These poles occur at the roots $\eta_{n,\hmu}$, $n=0,1,2,\ldots$, of the equation (recall that $\rA_0\equiv\Ai$)
\begin{align}
\label{eqn:ImpedanceRoots}
\hmu\re^{\ri\pi/3}\rA_0(\eta)+\rA_0'(\eta) = 0.
\end{align}
When $\hmu=\infty$ (Dirichlet) or $\hmu=0$ (Neumann) the roots $\eta_{n,\hmu}$ all lie on the negative real axis, being respectively the roots of $\rA_0=\Ai$ and $\rA_0'=\Ai'$, and the picture is as illustrated schematically in \F\ref{fig:gamma}. In the general case they asymptote to the negative real axis as $n\to\infty$, since $\rA_0'$ is the dominant term in \rf{eqn:ImpedanceRoots} for large $\eta$. 
\begin{align}
\hV(t,\hmu) = \frac{\re^{-\ri 2\pi/3}}{2\pi} \sum_{n=0}^\infty \frac{(\hmu^2+\re^{\ri\pi/3}\eta_{n,\hmu})\re^{\re^{-\ri\pi/6}t\eta_{n,\hmu}}}{\left(\hmu\rA_0'(\eta_{n,\hmu})+\re^{-\ri\pi/3}\eta_{n,\hmu}\rA_0(\eta_{n,\hmu})\right)^2}, \qquad \arg{t}\in(-\pi/3,2\pi/3).
\label{eqn:PekerisResidueImp}
\tagR{eqn:PekerisResidue}
\end{align}
\begin{align}
\hV(t,\hmu) \sim \frac{\re^{-\ri 2\pi/3}}{2\pi} \frac{(\hmu^2+\re^{\ri\pi/3}\eta_{0,\hmu})\re^{\re^{-\ri\pi/6}t\eta_{0,\hmu}}}{\left(\hmu\rA_0'(\eta_{0,\hmu})+\re^{-\ri\pi/3}\eta_{0,\hmu}\rA_0(\eta_{0,\hmu})\right)^2}, \qquad |t|\to\infty,\,\arg{t}\in(-\pi/3,2\pi/3).
\label{eqn:pPekerisAsympt1Imp}
\tagR{eqn:pPekerisAsympt1}
\end{align}
\begin{align}
\hV(t,\hmu)\sim -\left( \frac{t/2-\ri\hmu}{t/2+\ri\hmu} \right)\frac{\sqrt{-t}}{2\sqrt{\pi}}\re^{-\ri (t^3/12-\pi/4)} \left(1+ \ord{\frac{1}{|t|^3}} \right), \quad |t|\to\infty,\,\arg{t} \in(2\pi/3,5\pi/3).
\label{eqn:pPekerisAsympt2Imp}
\tagR{eqn:pPekerisAsympt2}
\end{align}
\begin{align}
A^s\sim- \frac{1}{\sqrt{3}}\left( \frac{\tau_-/2-\ri\mu/k}{\tau_-/2+\ri\mu/k}\right)\left(1-\frac{x}{\sqrt{x^2+3y}} \right)^{1/2}\re^{\ri k (4/27)\left(-x^3-(9/2)xy+(x^2+3y)^{3/2}\right)},
\label{eqn:AsMatchingIlluminatedImp}
\tagR{eqn:AsMatchingIlluminated}
\end{align}
Note that the factor $(\tau_-/2-\ri\mu/k)/(\tau_-/2+\ri\mu/k)$ is the inner limit of the usual reflection coefficient $(\cos{\theta}+\ri \mu/k)/(\cos{\theta}-\ri \mu/k)$ for reflection by an impedance boundary (where $\theta$ is the acute angle between the incident ray and the normal to the boundary). 
\begin{align}
g(\xi) = \sqrt{2\pi}\re^{\ri 3\pi/4}\re^{-\ri \xi^3/3}\, \hV(-\xi,\hmu).
\label{eqn:gDefImp}
\tagR{eqn:gDef}
\end{align}
\begin{align}
\tg(\xi) = \sqrt{2\pi}\re^{\ri 3\pi/4}\re^{-\ri \xi^3/3}\left(\hV(-\xi,\hmu)+\frac{\re^{\ri \xi^3/3}}{2\pi\ri\xi} \right)
&= \sqrt{2\pi}\re^{\ri 3\pi/4}\re^{-\ri \xi^3/3}\left(V(-\xi,\hmu)+\frac{\left(\re^{\ri \xi^3/3}-1\right)}{2\pi\ri\xi} \right).
\label{eqn:tgDefImp}
\tagR{eqn:tgDef}
\end{align}
\begin{align}
A\sim \frac{(\hmu^2+\re^{\ri\pi/3}\eta_{0,\hmu})\re^{-\ri k (xy/2+x^3/12)}\re^{\ri k^{1/3} \re^{\ri\pi/3}\eta_{0,\hmu} x}}{\left(\hmu \rA_0'(\eta_{0,\hmu})+\re^{-\ri\pi/3}\eta_{0,\hmu}\rA_0(\eta_{0,\hmu})\right)^2}\rA_0(\eta_{0,\hmu}+\re^{-\ri\pi/3}\hn),
\label{eqn:MatchingCreepingImp}
\tagR{eqn:MatchingCreeping}
\end{align}

\subsection{Neumann case}\label{sec:Neumann}
\begin{align}
\label{eqn:NeumannBC}
\pdone{\phi}{\bn}=0, \qquad \textrm{on }\partial D.
\tagN{eqn:DirichletBC}
\end{align}
\begin{align}
\label{eqn:FLNeuBC}
\pdone{A}{\hy}+\frac{\ri\hx}{2}A=0, \qquad \textrm{on } \hy+\frac{\hx^2}{4}=0.
\tagN{eqn:FLDirichletBC}
\end{align}
\begin{align}
\label{eqn:FLScatClassicalNeu}
\tagN{eqn:FLScatClassical}
A^s &= -\re^{-\ri(\hx \hy/2+ \hx^3/12)}\int_{-\infty}^\infty\re^{\ri \hx\sigma/2} \frac{\rA_0'(\sigma)}{ \rA_1'(\sigma)}\rA_1(\sigma-\hn) \,\rd\sigma.\\
A &= \re^{-\ri(\hx \hy/2+ \hx^3/12)}\int_{-\infty}^\infty\re^{\ri \hx\sigma/2}\left[\rA_0(\sigma-\hn)-\frac{\rA_0'(\sigma)}{ \rA_1'(\sigma)}\rA_1(\sigma-\hn)\right] \,\rd\sigma.
\label{eqn:FLTotClassicalNeu}
\tagN{eqn:FLTotClassical}
\end{align}
\begin{align}
\label{}
A^s &= \re^{-\ri(\hx \hy/2+ \hx^3/12)}
\left(\int_{l_1} \re^{\ri \hx\sigma/2}\rA_1(\sigma-\hn)\,\rd \sigma
+ \int_{l_2} \re^{\ri \hx\sigma/2}\frac{\rA_2'(\sigma)}{ \rA_1'(\sigma)}\rA_1(\sigma-\hn) \,\rd \sigma \right. \notag \\
&\hs{50} \left. - \int_{l_3} \re^{\ri \hx\sigma/2}\frac{\rA_0'(\sigma)}{ \rA_1'(\sigma)}\rA_1(\sigma-\hn) \,\rd \sigma\right)
\label{eqn:FLScatHornyNeu}
\tagN{eqn:FLScatHorny}
\end{align}
\begin{align}
\label{eqn:FLScatNewNeu}
A^s = \int_{\Gamma_1^l} \hq(t) \re^{\ri(-\hy t-\hx t^2/2 + t^3/3)}\,\rd t,
\tagN{eqn:FLScatNew}
\end{align}
\begin{align}
\hq(t)\equiv \hV(t,0) &= \frac{1}{2\pi}\left(\int_{l_2} \re^{\ri t\sigma}\,\rd \sigma
+ \int_{l_2} \re^{\ri t\sigma}\frac{\rA_2'(\sigma)}{\rA_1'(\sigma)} \,\rd \sigma
- \int_{l_3} \re^{\ri t\sigma}\frac{\rA_0'(\sigma)}{\rA_1'(\sigma)} \,\rd \sigma\right),\notag\\
& \hs{75} \arg{t}\in(-2\pi/3,\pi/3).
\label{eqn:PekerisRep1cNeu}
\tagN{eqn:PekerisRep1c}
\end{align}
\begin{align}
\hq(t) = \frac{1}{2\pi\ri t} + q(t),
\qquad t\neq 0,
\label{eqn:PekerisRep1bNeu}
\tagN{eqn:PekerisRep1b}
\end{align}
\begin{align}
q(t)=\frac{1}{2\pi}\left(\int_{l_2} \re^{\ri t\sigma}\frac{\rA_2'(\sigma)}{\rA_1'(\sigma)} \,\rd \sigma
- \int_{l_3} \re^{\ri t\sigma}\frac{\rA_0'(\sigma)}{\rA_1'(\sigma)} \,\rd \sigma\right).
\label{eqn:PekerisDefnNeu}
\tagN{eqn:PekerisDefn}
\end{align}
\begin{align}
A = \int_{\Gamma_1^r} \hq(t) \re^{\ri(-\hy t-\hx t^2/2 + t^3/3)}\,\rd t,
\label{eqn:FLTotNewNeu}
\tagN{eqn:FLTotNew}
\end{align}
\begin{align}
\hq(t)  = -\frac{1}{2\pi}\int_{-\infty}^\infty \re^{\ri t\sigma}\frac{\rA_0'(\sigma)}{\rA_1'(\sigma)} \,\rd \sigma,
\quad \im{t}<0.
\label{eqn:PekerisRep2Neu}
\tagN{eqn:PekerisRep2}
\end{align}
\begin{align}
\hq(t) = \frac{1}{4\pi^{2}t}\int_{L} \frac{\eta\re^{\re^{-\ri\pi/6}t\eta}}{\rA_0'(\eta)^2} \,\rd \eta, 
\label{eqn:pPekerisNeu}
\tagN{eqn:pPekeris}
\end{align}
\begin{align}
\hq(t)= -\frac{\re^{-\ri 2\pi/3}}{2\pi} \sum_{n=0}^\infty \frac{\re^{\re^{-\ri\pi/6}t\eta_{n,0}}}{\eta_{n,0}\rA_0(\eta_{n,0})^2}, \qquad \arg{t}\in(-\pi/3,2\pi/3).
\label{eqn:PekerisResidueNeu}
\tagN{eqn:PekerisResidue}
\end{align}
\begin{align}
\hq(t) \sim -\frac{\re^{-\ri 2\pi/3}}{2\pi} \frac{\re^{\re^{-\ri\pi/6}t\eta_{0,0}}}{\eta_{0,0}\rA_0(\eta_{0,0})^2}, \qquad |t|\to\infty,\,\arg{t}\in(-\pi/3,2\pi/3).
\label{eqn:pPekerisAsympt1Neu}
\tagN{eqn:pPekerisAsympt1}
\end{align}
\begin{align}
\hq(t)\sim -\frac{\sqrt{-t}}{2\sqrt{\pi}}\re^{-\ri (t^3/12-\pi/4)} \left(1+ \ord{\frac{1}{|t|^3}} \right), \quad |t|\to\infty,\,\arg{t} \in(2\pi/3,5\pi/3).
\label{eqn:pPekerisAsympt2Neu}
\tagN{eqn:pPekerisAsympt2}
\end{align}
\begin{align}
A^s\sim - \frac{1}{\sqrt{3}}\left(1-\frac{x}{\sqrt{x^2+3y}} \right)^{1/2}\re^{\ri k (4/27)\left(-x^3-(9/2)xy+(x^2+3y)^{3/2}\right)},
\label{eqn:AsMatchingIlluminatedNeu}
\tagN{eqn:AsMatchingIlluminated}
\end{align}
\begin{align}
g(\xi) = \sqrt{2\pi}\re^{\ri 3\pi/4}\re^{-\ri \xi^3/3}\, \hq(-\xi).
\label{eqn:gDefNeu}
\tagN{eqn:gDef}
\end{align}
\begin{align}
\tg(\xi) = \sqrt{2\pi}\re^{\ri 3\pi/4}\re^{-\ri \xi^3/3}\left(\hq(-\xi)+\frac{\re^{\ri \xi^3/3}}{2\pi\ri\xi} \right)
&= \sqrt{2\pi}\re^{\ri 3\pi/4}\re^{-\ri \xi^3/3}\left(q(-\xi)+\frac{\left(\re^{\ri \xi^3/3}-1\right)}{2\pi\ri\xi} \right).
\label{eqn:tgDefNeu}
\tagN{eqn:tgDef}
\end{align}
\begin{align}
A\sim- \frac{\re^{-\ri k (xy/2+x^3/12)}\re^{\ri k^{1/3} \re^{\ri\pi/3}\eta_{0,0} x}}{\eta_{0,0}\rA_0(\eta_{0,0})^2}\rA_0(\eta_{0,0}+\re^{-\ri\pi/3}\hn),
\label{eqn:MatchingCreepingNeu}
\tagN{eqn:MatchingCreeping}
\end{align}

\section{Airy function notation}
\label{app:Airy}

\begin{figure}[t]
\colorlet{lightgray}{black!15}
\begin{center}
\subfigure[(a)]{
\begin{tikzpicture}[line cap=round,line join=round,>=triangle 45,x=1.0cm,y=1.0cm, scale=2]
\def\xm{0.75};
\def\ym{0.75};
\pgfmathsetmacro{\rad}{2*sqrt(\xm^2+\ym^2)}
\def\theta{60};
\pgfmathsetmacro{\costheta}{cos(\theta)};
\pgfmathsetmacro{\sintheta}{sin(\theta)};
\pgfmathsetmacro{\costwotheta}{cos(2*\theta)};
\pgfmathsetmacro{\sintwotheta}{sin(2*\theta)};
\pgfmathsetmacro{\costhreetheta}{cos(3*\theta)};
\pgfmathsetmacro{\sinthreetheta}{sin(3*\theta)};
\pgfmathsetmacro{\cosfourtheta}{cos(4*\theta)};
\pgfmathsetmacro{\sinfourtheta}{sin(4*\theta)};
\pgfmathsetmacro{\cosfivetheta}{cos(5*\theta)};
\pgfmathsetmacro{\sinfivetheta}{sin(5*\theta)};
\pgfmathsetmacro{\cossixtheta}{cos(6*\theta)};
\pgfmathsetmacro{\sinsixtheta}{sin(6*\theta)};
\begin{scope}
\clip (-\xm,-\ym) rectangle (\xm,\ym);
	\filldraw [lightgray] (0,0) -- (\rad,0) -- (\rad*\costheta,\rad*\sintheta) ;
	\filldraw [lightgray] (0,0) -- (\rad*\costwotheta,\rad*\sintwotheta) -- (\rad*\costhreetheta,\rad*\sinthreetheta) ;
	\filldraw [lightgray] (0,0) -- (\rad*\cosfourtheta,\rad*\sinfourtheta) -- (\rad*\cosfivetheta,\rad*\sinfivetheta) ;
\begin{scope}[scale=0.125,>=latex]
\draw [samples=50,domain=-0.99:0.99,rotate around={90:(0,0)},xshift=0cm,yshift=0cm] plot ({1*(1+(\x)^2)/(1-(\x)^2)},{1.73*2*(\x)/(1-(\x)^2)});
\draw [->] (1.4,1.27) -- (1.5,1.3);
\begin{scope}[rotate around={120:(0,0)}]
\draw [samples=50,domain=-0.99:0.99,rotate around={90:(0,0)},xshift=0cm,yshift=0cm] plot ({1*(1+(\x)^2)/(1-(\x)^2)},{1.73*2*(\x)/(1-(\x)^2)});
\draw [->] (1.4,1.27) -- (1.5,1.3);
\end{scope}
\begin{scope}[rotate around={240:(0,0)}]
\draw [samples=50,domain=-0.99:0.99,rotate around={90:(0,0)},xshift=0cm,yshift=0cm] plot ({1*(1+(\x)^2)/(1-(\x)^2)},{1.73*2*(\x)/(1-(\x)^2)});
\draw [->] (1.4,1.27) -- (1.5,1.3);
\end{scope}
\end{scope}
\end{scope}
\draw (\xm+0.15,-0.1) node {$\real{ t}$};
\draw [->] (-\xm,0) -- (\xm,0);
\draw [->] (0,-\ym) -- (0,\ym);
\draw (0.25,\ym-0.05) node {$\im{ t}$};
\draw (0.3*\xm,0.42*\ym) node {$\Gamma_{0}$};
\draw (-0.3*\xm,-0.2*\ym) node {$\Gamma_{1}$};
\draw (0.35*\xm,-0.3*\ym) node {$\Gamma_{2}$};
\end{tikzpicture} 
}
\subfigure[(b)]{
\begin{tikzpicture}[line cap=round,line join=round,>=triangle 45,x=1.0cm,y=1.0cm, scale=2]
\def\xm{0.75};
\def\ym{0.75};
\pgfmathsetmacro{\rad}{2*sqrt(\xm^2+\ym^2)}
\def\theta{60};
\pgfmathsetmacro{\costheta}{cos(\theta)};
\pgfmathsetmacro{\sintheta}{sin(\theta)};
\pgfmathsetmacro{\costwotheta}{cos(2*\theta)};
\pgfmathsetmacro{\sintwotheta}{sin(2*\theta)};
\pgfmathsetmacro{\costhreetheta}{cos(3*\theta)};
\pgfmathsetmacro{\sinthreetheta}{sin(3*\theta)};
\pgfmathsetmacro{\cosfourtheta}{cos(4*\theta)};
\pgfmathsetmacro{\sinfourtheta}{sin(4*\theta)};
\pgfmathsetmacro{\cosfivetheta}{cos(5*\theta)};
\pgfmathsetmacro{\sinfivetheta}{sin(5*\theta)};
\pgfmathsetmacro{\cossixtheta}{cos(6*\theta)};
\pgfmathsetmacro{\sinsixtheta}{sin(6*\theta)};
\def\sc{0.15}
\def\nn{0.4}
\fill (\sc*2.34*\costhreetheta,\sc*2.34*\sinthreetheta) circle (1pt);
\fill (\sc*4.09*\costhreetheta,\sc*4.09*\sinthreetheta) circle (1pt);
\fill (\sc*5.52*\costhreetheta,\sc*5.52*\sinthreetheta) circle (1pt);
\fill (\sc*6.79*\costhreetheta,\sc*6.79*\sinthreetheta) circle (1pt);
\draw (\xm+0.2,-0.1) node {$\real{z}$};
\draw [->] (-\xm,0) -- (\xm,0);
\draw [->] (0,-\ym) -- (0,\ym);
\draw (0.3,\ym-0.05) node {$\im{z}$};
\draw [<->,line width=0.5pt] (0.5*\costheta,-0.5*\sintheta) arc[x radius=0.5cm, y radius =.5cm, start angle=-60, end angle=60];
\draw [<->,line width=0.5pt] (-0.5,0) arc[x radius=0.5cm, y radius =.5cm, start angle=-180, end angle=-60];
\draw [<->,line width=0.5pt] (-0.5,0) arc[x radius=0.5cm, y radius =.5cm, start angle=180, end angle=60];
\draw (0.7,0.25) node {decay};
\draw (-0.55,-0.5) node {growth};
\draw (-0.55,0.5) node {growth};
\draw (-0.80,0.1) node {zeros};
\draw [dotted] (0,0) -- (\xm*\costheta,\xm*\sintheta);
\draw [dotted] (0,0) -- (\xm*\costheta,-\xm*\sintheta);
\end{tikzpicture} 
}
\caption{(a) The contours $\Gamma_j$, $j=0,1,2$, and the sectors (shaded) in which $\re^{\ri t^3}$ decays exponentially as $|t|\to\infty$. (b) Schematic showing the large-argument asymptotic behaviour (exponential growth or decay) of the Airy function $\rA_0(z)=\Ai(z)$, and the zeros on the negative real axis. The same qualitative behaviour is exhibited by the derivative $\rA_0(z)'=\Ai'(z)$.}
\label{Sectors}
\end{center}
\end{figure}
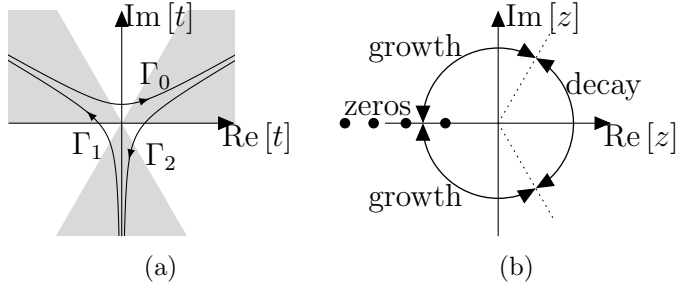 

We consider the following solutions of the Airy equation $A''-zA=0$:
\begin{align}
\label{eqn:AjDef}
\rA_j(z) := \frac{1}{2\pi}\int_{\Gamma_j}\re^{\ri (t z+t^3/3)}\,\rd t, \qquad z\in \C, \,\,j=0,1,2,
\end{align}
where, for each $j=0,1,2$, the integration contour $\Gamma_j$ runs from infinity at $\arg{t}=(4j+5)\pi/6$ to infinity at $\arg{t}=(4j+1)\pi/6$, as illustrated in Fig.~\ref{Sectors}(a). The functions $\rA_j$ are related to other common Airy function notation as follows:
\begin{align*}
\label{}
\rA_0(z) &= \Ai(z) = v(z)/\sqrt{\pi},\\
\rA_1(z) &= \re^{\ri 2\pi/3}\Ai(\re^{\ri 2\pi/3} z) = \frac{\ri}{2\sqrt{\pi}}w_1(z),\\
\rA_2(z) &= \re^{-\ri 2\pi/3}\Ai(\re^{-\ri 2\pi/3} z) = -\frac{\ri}{2\sqrt{\pi}}w_2(z),
\end{align*}
where $\Ai$ is the usual Airy function and $v,w_1,w_2$ are the Airy functions commonly found in the Russian literature (cf.\ e.g.\ \cite{BaKi:79,BaBu:91}). The functions $\rA_j$ satisfy the connection formula
\begin{align}
\label{eqn:Connection}
\sum_{j=0}^2 \rA_j(z) = 0, \qquad z\in\C,
\end{align}
and the Wronskian relation\footnote{We remark that the Wronskian formulas in \cite[p.~405]{BaBu:91} are incorrect: instead of $\{2,\ri,\ri,1\}$ they should read $\{2\ri,-1,-1,-1\}$.} 
 (cf.\ \cite[\S9.2]{DLMF})
\begin{align}
\label{eqn:AiryWronskian}
W(\rA_j,\rA_{j+1})=\rA_{j+1}'\rA_j-\rA_{j+1}\rA_j' = \frac{\ri}{2\pi}, \qquad j=0,1,2, \qquad \rA_3:=\rA_0.
\end{align}
The following large argument ($|z| \to \infty$) behaviour is well known (see e.g.\ \cite[\S9.7]{DLMF}):
\begin{align}
\Ai(z) & \sim \frac{\re^{-(2/3) z^{3/2}}}{2\sqrt{\pi}z^{1/4}}\left(1+\ord{\frac{1}{|z|^{3/2}}}\right), 
\qquad |\arg{z}|\leq \pi-\delta, \label{eqn:AiryAsympt1} \\
\Ai(-z) & \sim \frac{\sin((2/3) z^{3/2}+\pi/4)}{\sqrt{\pi}z^{1/4}}\left(1+\ord{\frac{1}{|z|^{3}}}\right)\notag\\
& \quad - \frac{5\cos((2/3) z^{3/2}+\pi/4)}{48\sqrt{\pi}z^{7/4}}\left(1+\ord{\frac{1}{|z|^{3}}}\right), 
\qquad |\arg{z}|\leq 2\pi/3-\delta, \label{eqn:AiryAsympt2} \\
\Ai'(z) & \sim -\frac{z^{1/4}\re^{-(2/3) z^{3/2}}}{2\sqrt{\pi}}\left(1+\ord{\frac{1}{|z|^{3/2}}}\right), 
\qquad |\arg{z}|\leq \pi-\delta, \label{eqn:AiryPrimeAsympt1}\\
\Ai'(-z) & \sim \frac{-z^{1/4}\sin((2/3) z^{3/2}+\pi/4)}{\sqrt{\pi}}\left(1+\ord{\frac{1}{|z|^{3}}}\right)\notag\\
& \quad + \frac{7\cos((2/3) z^{3/2}+\pi/4)}{48\sqrt{\pi}z^{5/4}}\left(1+\ord{\frac{1}{|z|^{3}}}\right), 
\qquad |\arg{z}|\leq 2\pi/3-\delta, \label{eqn:AiryPrimeAsympt2} 
\end{align}
where $\delta$ is an arbitrary positive constant and the approximations hold uniformly in the ranges of $\arg{z}$ stated, with the principle branches of all multi-valued functions assumed throughout. The qualitative asymptotic behaviour is illustrated schematically in \F\ref{Sectors}(b).

Finally, we note the following well-known Fourier transform relation \cite[eqn.~6.1.22]{BaKi:79}
\begin{align}
\label{eqn:AiryRelation3}
\int_{-\infty}^\infty \re^{-\ri z t} \Ai(z) \,\rd z = \re^{\ri t^3/3}, \qquad t\in \R,
\end{align}
where the integral is understood in an improper sense, and also the following integral identity, which follows from a simple change of variable in \rf{eqn:AjDef} and holds for all $\hx,\hy,\sigma \in \C$:
\begin{align}
\label{eqn:AiryIdentity}
\rA_j(\sigma-(\hy+\frac{\hx^2}{4})) =\frac{1}{2\pi}\re^{-\ri \hx\sigma/2} \re^{\ri(\hx\hy/2+\hx^3/12)}\int_{\Gamma_j}\re^{\ri \sigma t}\re^{\ri(-\hy t-\hx t^2/2 + t^3/3)} \, \rd t, \quad j=0,1,2.
\end{align}
\fi
\section{Matching the classical solution to the penumbra field}
\label{app:JRO}
In this appendix [supplied by J.\ Ockendon] we sketch how the representation \rf{eqn:FLScatClassical} can be matched to the penumbra field using the method of stationary phase. Writing \rf{eqn:FLScatClassical} in terms of the variables $x=k^{-1/3}\hx$, $\ty=k^{-1/3}\hy$ gives
\begin{align}
\label{}
A^s &= -\re^{-\ri(k^{2/3}x \ty/2+k x^3/12)}\int_{-\infty}^\infty\re^{\ri k^{1/3} x\sigma/2} \frac{\rA_0(\sigma)}{\rA_1(\sigma)}\rA_1(\sigma-k^{2/3}n) \,\rd\sigma,
\label{eqn:FLScatClassical2}
\end{align}
where $n:=x^2/4+k^{-1/3}\ty$. To determine the behaviour of \rf{eqn:FLScatClassical2} as $k\to\infty$ we split the integral into three parts (where $\Sigma,\Sigma'\gg 1$ are as yet unspecified):
\begin{align*}
\label{}
\int_{-\infty}^\infty \,\rd\sigma = \int_{-\infty}^{-\Sigma} \,\rd\sigma + \int_{-\Sigma}^{\Sigma'} \,\rd\sigma + \int_{\Sigma'}^\infty  \,\rd\sigma
=I_1 + I_2 + I_3.
\end{align*}
$I_3$ is exponentially small as $\Sigma'\to\infty$, and can be neglected. In $I_2$ and $I_1$ we replace $\rA_1(\sigma-k^{2/3}n)$ by its large negative argument approximation 
\begin{align}
\label{eqn:Airy1exp}
\rA_1(\sigma-k^{2/3}n) \sim \frac{\re^{\ri 3\pi/4}\re^{(2\ri/3)(-\sigma+k^{2/3}n)^{3/2}}}{2\sqrt{\pi}(-\sigma+k^{2/3}n)^{1/4}}, \qquad \sigma-k^{2/3}n \to-\infty,
\end{align}
which applies in $I_2$ provided $\Sigma'\ll k^{2/3}$. By further approximating \rf{eqn:Airy1exp} one can show that, provided $\Sigma \ll k^{1/6}$ (so that terms of order $\sigma^2$ and above can be neglected when expanding the phase), $I_2$ makes a contribution to $A^s$ of
\begin{align}
\label{eqn:I2Def}
\frac{\sqrt{2\pi}\re^{3\ri \pi/4}}{k^{1/6}x^{1/2}}\re^{\ri k^{1/3} \ty^2/(2x) }\re^{-\ri \ty^3/(3x^3)} I_\Sigma\left(-\frac{\ty}{x}\right),
\end{align}
where
\begin{align*}
\label{}
I_\Sigma(t):= -\frac{1}{2\pi}\int_{-\Sigma}^\infty \re^{\ri t\sigma}\frac{\rA_0(\sigma)}{\rA_1(\sigma)} \,\rd \sigma.
\end{align*}
(Taking the upper integration limit as $\infty$ introduces only exponentially small errors.) The integral $I_\Sigma(t)$ closely resembles the representation \rf{eqn:PekerisRep2} for the Pekeris caret function. In fact, one can check using contour integration that 
\begin{align}
\label{eqn:ISigmaDef}
I_\Sigma(t) = \hp(t) - \frac{\re^{-\ri t \Sigma}}{2\pi\ri t} + R_\Sigma(t), 
\end{align}
where the remainder $R_\Sigma(t)$ is $\ord{\Sigma^{-1/2}}$ as $\Sigma\to\infty$, for $t\geq -\Sigma$. In $I_1$ we can approximate $\rA_0(\sigma)/\rA_1(\sigma)\sim -1$, plus a rapidly oscillating term which we ignore since it contributes a higher order correction to the integral. Rescaling $\sigma=k^{1/3}\tsigma$ and writing, with $1\ll\tSigma\ll k^{1/12}$,
\begin{align*}
\label{}
I_1 = \int_{-\infty}^\Sigma \,\rd\sigma = k^{1/3}\int_{-\infty}^{-\tSigma} \,\rd\tsigma + k^{1/3}\int_{-\tSigma}^{-k^{-1/3}\Sigma} \,\rd\tsigma =I_{11}+ I_{12}, 
\end{align*}
we neglect the contribution from $I_{11}$ and expand the phase in $I_{12}$ to obtain a contribution to $A^s$ of
\begin{align*}
\label{}
\frac{k^{1/6}\re^{3\ri \pi/4}}{\sqrt{2\pi}x^{1/2}}
\int_{-\tSigma}^{-k^{-1/3}\Sigma} \re^{\ri\left( k^{1/3} (\ty-\tsigma)^2/(2x) -(\ty-\tsigma)^3/(3x^3)\right)} \,\rd\tsigma.
\end{align*}
When $\ty<0$ and $|\ty|>k^{-1/3} \Sigma$ there is a stationary phase point at $\tsigma=\ty$ which makes a contribution of $-1$ to $A^s$, cancelling the incident field. The endpoint contribution from $\tsigma=k^{-1/3} \Sigma$ makes a contribution to $A^s$ of
\begin{align*}
\label{}
\frac{\ri\re^{3\ri \pi/4}x^{1/2}}{\sqrt{2\pi}k^{1/6}\ty}
\re^{\ri k^{1/3} \ty^2/(2x) }\re^{-\ri \ty^3/(3x^3)}  \re^{\ri \Sigma\ty/x},
\end{align*}
which (as it must) cancels the $\Sigma$-dependent term in \rf{eqn:I2Def} arising from the second term on the right-hand side of \rf{eqn:ISigmaDef}. We have thus recovered the behaviour of the field in the transition regions V$_\textrm{upper}$ and V$_\textrm{lower}$. For $\ty=\ord{k^{-1/6}}$ (i.e.\ in the Fresnel region IV), by rescaling $\ty=k^{-1/6}\check y$ and $\tsigma=k^{-1/6}\check \sigma$, we see that $I_{12}$ gives a contribution
\begin{align*}
\label{}
\frac{\re^{3\ri \pi/4}}{\sqrt{2\pi}x^{1/2}}
\int_{-k^{1/6}\tSigma}^{-k^{-1/6}\Sigma} \re^{\ri k^{1/3} (\check y-\check \sigma)^2/(2x)} \,\rd\check \sigma 
\sim \Fr\left(-\frac{\check y}{\sqrt{2x}}\right),
\end{align*}
so that the emergence of the Fresnel integral in this region is also verified.

%
\end{document}

%% file: macros.tex
\newcommand{\done}[2]{\dfrac{d {#1}}{d {#2}}}
\newcommand{\donet}[2]{\frac{d {#1}}{d {#2}}}
\newcommand{\pdone}[2]{\dfrac{\partial {#1}}{\partial {#2}}}
\newcommand{\pdonet}[2]{\frac{\partial {#1}}{\partial {#2}}}
\newcommand{\pdonetext}[2]{\partial {#1}/\partial {#2}}
\newcommand{\pdtwo}[2]{\dfrac{\partial^2 {#1}}{\partial {#2}^2}}
\newcommand{\pdtwot}[2]{\frac{\partial^2 {#1}}{\partial {#2}^2}}
\newcommand{\pdtwomix}[3]{\dfrac{\partial^2 {#1}}{\partial {#2}\partial {#3}}}
\newcommand{\pdtwomixt}[3]{\frac{\partial^2 {#1}}{\partial {#2}\partial {#3}}}
\newcommand{\bs}[1]{\mathbf{#1}}
\newcommand{\bx}{\mathbf{x}}
\newcommand{\by}{\mathbf{y}}
\newcommand{\bd}{\mathbf{d}} 
\newcommand{\bn}{\mathbf{n}} 
\newcommand{\bP}{\mathbf{P}} 
\newcommand{\bp}{\mathbf{p}} 
\newcommand{\ol}[1]{\overline{#1}}
\newcommand{\rf}[1]{(\ref{#1})}
\newcommand{\xt}{\mathbf{x},t}
\newcommand{\hs}[1]{\hspace{#1mm}}
\newcommand{\vs}[1]{\vspace{#1mm}}
\newcommand{\eps}{\varepsilon}
\newcommand{\ord}[1]{\mathcal{O}\left(#1\right)} 
\newcommand{\oord}[1]{o\left(#1\right)}
\newcommand{\Ord}[1]{\Theta\left(#1\right)}
\newcommand{\PhiF}{\Phi_{\rm freq}}
\newcommand{\real}[1]{{\rm Re}\left[#1\right]} 
\newcommand{\im}[1]{{\rm Im}\left[#1\right]}
\newcommand{\hsnorm}[1]{||#1||_{H^{s}(\bs{R})}}
\newcommand{\hnorm}[1]{||#1||_{\tilde{H}^{-1/2}((0,1))}}
\newcommand{\norm}[2]{\left\|#1\right\|_{#2}}
\newcommand{\normt}[2]{\|#1\|_{#2}}
\newcommand{\on}[1]{\Vert{#1} \Vert_{1}}
\newcommand{\tn}[1]{\Vert{#1} \Vert_{2}}
\newcommand{\ts}{\tilde{s}}
\newcommand{\darg}[1]{\left|{\rm arg}\left[ #1 \right]\right|}
\newcommand{\bnabla}{\boldsymbol{\nabla}}
\newcommand{\dive}{\boldsymbol{\nabla}\cdot}
\newcommand{\curl}{\boldsymbol{\nabla}\times}
\newcommand{\Phixy}{\Phi(\bx,\by)}
\newcommand{\PhiOxy}{\Phi_0(\bx,\by)}
\newcommand{\dxPhixy}{\pdone{\Phi}{n(\bx)}(\bx,\by)}
\newcommand{\dyPhixy}{\pdone{\Phi}{n(\by)}(\bx,\by)}
\newcommand{\dxPhiOxy}{\pdone{\Phi_0}{n(\bx)}(\bx,\by)}
\newcommand{\dyPhiOxy}{\pdone{\Phi_0}{n(\by)}(\bx,\by)}

\newcommand{\rd}{\mathrm{d}}
\newcommand{\R}{\mathbb{R}}
\newcommand{\N}{\mathbb{N}}
\newcommand{\Z}{\mathbb{Z}}
\newcommand{\C}{\mathbb{C}}
\newcommand{\K}{{\mathbb{K}}}
\newcommand{\ri}{{\mathrm{i}}}
\newcommand{\re}{{\mathrm{e}}} 

\newcommand{\cA}{\mathcal{A}}
\newcommand{\cC}{\mathcal{C}}
\newcommand{\cS}{\mathcal{S}}
\newcommand{\cD}{\mathcal{D}}
\newcommand{\cone}{{c_{j}^\pm}}
\newcommand{\ctwo}{{c_{2,j}^\pm}}
\newcommand{\cthree}{{c_{3,j}^\pm}}

\newtheorem{thm}{Theorem}[section]
\newtheorem{lem}[thm]{Lemma}
\newtheorem{defn}[thm]{Definition}
\newtheorem{prop}[thm]{Proposition}
\newtheorem{cor}[thm]{Corollary}
\newtheorem{rem}[thm]{Remark}
\newtheorem{conj}[thm]{Conjecture}
\newtheorem{ass}[thm]{Assumption}
\newtheorem{example}[thm]{Example} 

%% file: Hewett_FLP_v1_preprint.bbl
\begin{thebibliography}{10}
\expandafter\ifx\csname url\endcsname\relax
  \def\url#1{\texttt{#1}}\fi
\expandafter\ifx\csname urlprefix\endcsname\relax\def\urlprefix{URL }\fi
\expandafter\ifx\csname href\endcsname\relax
  \def\href#1#2{#2} \def\path#1{#1}\fi

\bibitem{Fo:46}
V.~A. Fock, The field of a plane wave near the surface of a conducting body, J.
  Phys. USSR 10 (1945) 399--409.

\bibitem{LeFo:45}
M.~S. Leontovich, V.~A. Fock, Solution of the problem of propagation of
  electromagnetic waves along the earth's surface by the method of the
  parabolic equation, J. Phys. USSR 10 (1946) 13--24.

\bibitem{Pe:47}
C.~L. Pekeris, The field of a microwave dipole antenna in the vicinity of the
  horizon, J. Appl. Phys. 18~(7) (1947) 667--680.

\bibitem{Fo:65}
V.~A. Fock, Electromagnetic Diffraction and Propagation Problems, Pergamon,
  Oxford, 1965.

\bibitem{BaKi:79}
V.~M. Babich, N.~Y. Kirpichnikova, The Boundary-layer Method in Diffraction
  Problems, Springer, Berlin, 1979.

\bibitem{BaBu:91}
V.~M. Babich, V.~S. Buldyrev, Short-Wavelength Diffraction Theory, Springer,
  Berlin, 1991.

\bibitem{TeChKiOcSmZa:00}
R.~H. Tew, S.~J. Chapman, J.~R. King, J.~R. Ockendon, B.~J. Smith,
  I.~Zafarullah, Scalar wave diffraction by tangent rays, Wave Motion 32 (2000)
  363--380.

\bibitem{Br:66}
W.~P. Brown~Jr, On the asymptotic behavior of electromagnetic fields scattered
  from convex cylinders near grazing incidence, J. Math. Anal. Appl. 15~(2)
  (1966) 355--385.

\bibitem{Lu:67}
D.~Ludwig, Uniform asymptotic expansion of the field scattered by a convex
  object at high frequencies, Comm. Pure Appl. Math. 20~(1) (1967) 103--138.

\bibitem{He:68}
A.~J. Hermans, High-frequency scattering by a convex smooth object, Ph.D.
  thesis, TU Delft (1968).

\bibitem{MeTa:86}
R.~B. Melrose, M.~E. Taylor, The radiation pattern of a diffracted wave near
  the shadow boundary, Commun. Part. Diff. Eq. 11~(6) (1986) 599--672.

\bibitem{BuLy:87}
V.~S. Buldyrev, M.~A. Lyalinov, Uniform and local asymptotic behavior of the
  penumbral wave field for diffraction of short waves on a smooth convex
  contour, J. Sov. Math. (now J. Math. Sci.) 38~(1) (1987) 1579--1584.

\bibitem{Lo:59}
N.~A. Logan, General Research in Diffraction Theory, Vol.~1, Lockheed Missiles
  and Space Division Technical Report, Lockheed Aircraft Corporation, 1959,
  available from \texttt{http://www.dtic.mil/get-tr-doc/pdf?AD=AD0241228}.

\bibitem{James}
G.~L. James, Geometrical Theory of Diffraction for Electromagnetic Waves,
  Institution of Electrical Engineers, Peter Peregrinus Ltd, 1986.

\bibitem{Pe:87}
L.~W. Pearson, A scheme for automatic computation of {F}ock-type integrals,
  IEEE Trans. Antennas Propag. 35~(10) (1987) 1111--1118.

\bibitem{Sa:09}
S.-E. Sandstr{\"o}m, Computation of the {F}ock scattering functions, in:
  Proceedings of the 3rd Conference on Mathematical Modeling of Wave Phenomena,
  20th Nordic Conference on Radio Science and Communications, Vol. 1106, AIP
  Publishing, 2009, pp. 104--109.

\bibitem{MeTa:85}
R.~B. Melrose, M.~E. Taylor, Near peak scattering and the corrected {K}irchhoff
  approximation for a convex obstacle, Adv. Math. 55~(3) (1985) 242--315.

\bibitem{SmithThesis}
B.~J. Smith, A complex ray approach to the acoustics of fluid-loaded
  structures, Ph.D. thesis, University of Nottingham (1995).

\bibitem{CoatsThesis}
J.~Coats, High frequency asymptotics of antenna/structure interactions, Ph.D.
  thesis, University of Oxford (2002).

\bibitem{FozardThesis}
J.~A. Fozard, Diffraction and scattering of high frequency waves, Ph.D. thesis,
  University of Oxford (2005).

\bibitem{EnKiTe:98}
J.~C. Engineer, J.~R. King, R.~H. Tew, Diffraction by slender bodies, Eur. J.
  Appl. Math 9 (1998) 129--158.

\bibitem{OckTew:12}
J.~R. Ockendon, R.~H. Tew, Thin-layer solutions of the {H}elmholtz and related
  equations, SIAM Rev. 54(1) (2012) 3--51.

\bibitem{DLMF}
{NIST} {D}igital {L}ibrary of {M}athematical {F}unctions,
  http://dlmf.nist.gov/, Release 1.0.7 of 2014-03-21.

\bibitem{Pop:79}
M.~M. Popov, The problem of whispering gallery waves in a neighbourhood of a
  simple zero of the effective curvature of the boundary, J. Sov. Math. (now J.
  Math. Sci.) 11 (1979) 791--797.

\bibitem{BabSmy:86}
V.~M. Babich, V.~P. Smyshlyaev, Scattering problems for the {S}chr{\"o}dinger
  equation in the case of a potential linear in time and coordinate {I}:
  {A}symptotics in the shadow zone, J. Sov. Math. (now J. Math. Sci.) 32 (1986)
  103--111.

\bibitem{Kaz:03}
A.~Y. Kazakov, Special function related to the concave-convex boundary problem
  of the diffraction theory, J. Phys. A: Math. Gen. 36 (2003) 4127--4141.

\end{thebibliography}
